\documentclass{article} %

\usepackage[
  margin=3cm,
  includefoot,
  footskip=30pt,
]{geometry}

\usepackage[utf8]{inputenc} %
\usepackage[T1]{fontenc}    %
\usepackage[pagebackref]{hyperref}
\usepackage{url}            %
\usepackage{booktabs}       %
\usepackage{amsfonts}       %
\usepackage{nicefrac}       %
\usepackage{microtype}      %
\usepackage{xcolor}         %

\definecolor{darkgreen}{rgb}{0.0,0.5,0.0}

\renewcommand*{\backrefalt}[4]{%
    \ifcase #1 \footnotesize{(Not cited.)}%
    \or        \footnotesize{(Cited on page~#2)}%
    \else      \footnotesize{(Cited on pages~#2)}%
    \fi}

\usepackage{amsfonts}
\usepackage{amsmath}
\usepackage{amsthm}
\usepackage{amssymb}
\usepackage{dsfont}
\usepackage{comment}
\usepackage{multirow}
\usepackage{multicol}
\usepackage{threeparttable}

\usepackage{xcolor}
\usepackage{color}
\usepackage{graphicx}
\usepackage{pifont}%
\newcommand{\algname}[1]{{\sf#1}\xspace}
\newcommand{\eqdef}{\coloneqq}

\usepackage{verbatim}
\usepackage{xspace} %

\usepackage{enumerate}
\usepackage{enumitem}

\usepackage{subcaption}

\providecommand{\algname}[1]{{\sf#1}\xspace}

\providecommand{\norm}[1]{\left\lVert#1\right\rVert}

\providecommand{\makecellnew}[1]{{\renewcommand{\arraystretch}{0.8}\begin{tabular}{c} #1 \end{tabular}}}

  \providecommand{\R}{\mathbb{R}} %

\providecommand{\EbNSq}[1]{\Eb{\norm{#1}^2}}
\providecommand{\EbInp}[2]{\Eb{\inp{#1}{#2}}}

  \DeclareMathOperator{\E}{{\mathbb E}}
  \providecommand{\wtilde}[1]{\widetilde{#1}}
  \providecommand{\Eb}[1]{{\mathbb E}\left[#1\right] }       %
  \providecommand{\EEb}[2]{{\mathbb E}_{#1}\left[#2\right] } %

  \DeclareMathOperator*{\argmin}{arg\,min}

  \providecommand{\0}{\mathbf{0}}
  
  \renewcommand{\aa}{\mathbf{a}}
  
  \providecommand{\cc}{\mathbf{c}}
  
  \providecommand{\ee}{\mathbf{e}}

  \renewcommand{\gg}{\mathbf{g}}

  \providecommand{\vv}{\mathbf{v}}
  
  \providecommand{\xx}{\mathbf{x}}
  \providecommand{\yy}{\mathbf{y}}

  \providecommand{\cC}{\mathcal{C}}
  \providecommand{\cD}{\mathcal{D}}

  \providecommand{\cO}{\mathcal{O}}

  \usepackage{bm}

\RequirePackage[colorinlistoftodos,bordercolor=orange,backgroundcolor=orange!20,linecolor=orange,textsize=scriptsize]{todonotes}

\providecommand{\mycomment}[3]{%
  \ifthenelse{\boolean{commentsOn}}{%
    \todo[caption={},color=#3!20,inline]{\textbf{#1: }#2}%
  }{}%
}
\providecommand{\myinlinecomment}[3]{%
  \ifthenelse{\boolean{commentsOn}}{%
    {\color{#1}#2: #3}%
  }{}%
}

\newcommand\commenter[2]{%
  \expandafter\newcommand\csname i#1\endcsname[1]{\myinlinecomment{#2}{#1}{##1}}
  \expandafter\newcommand\csname #1\endcsname[1]{\mycomment{#1}{##1}{#2}}
}

\usepackage{url}

\renewcommand{\epsilon}{\varepsilon}

\usepackage{algorithm}
\usepackage[noend]{algpseudocode}

\errorcontextlines\maxdimen

\makeatletter
    \newcommand*{\algrule}[1][\algorithmicindent]{\makebox[#1][l]{\hspace*{.5em}\thealgruleextra\vrule height \thealgruleheight depth \thealgruledepth}}%
\newcommand*{\thealgruleextra}{}
\newcommand*{\thealgruleheight}{.75\baselineskip}
\newcommand*{\thealgruledepth}{.25\baselineskip}

\newcount\ALG@printindent@tempcnta
\def\ALG@printindent{%
    \ifnum \theALG@nested>0%
        \ifx\ALG@text\ALG@x@notext%
        \else
            \unskip
            \addvspace{-1pt}%
            \ALG@printindent@tempcnta=1
            \loop
                \algrule[\csname ALG@ind@\the\ALG@printindent@tempcnta\endcsname]%
                \advance \ALG@printindent@tempcnta 1
            \ifnum \ALG@printindent@tempcnta<\numexpr\theALG@nested+1\relax%
            \repeat
        \fi
    \fi
    }%
\usepackage{etoolbox}
\patchcmd{\ALG@doentity}{\noindent\hskip\ALG@tlm}{\ALG@printindent}{}{\errmessage{failed to patch}}
\makeatother

\newbox\statebox
\newcommand{\myState}[1]{%
    \setbox\statebox=\vbox{#1}%
    \edef\thealgruleheight{\dimexpr \the\ht\statebox+1pt\relax}%
    \edef\thealgruledepth{\dimexpr \the\dp\statebox+1pt\relax}%
    \ifdim\thealgruleheight<.75\baselineskip
        \def\thealgruleheight{\dimexpr .75\baselineskip+1pt\relax}%
    \fi
    \ifdim\thealgruledepth<.25\baselineskip
        \def\thealgruledepth{\dimexpr .25\baselineskip+1pt\relax}%
    \fi
    \State #1%
    \def\thealgruleheight{\dimexpr .75\baselineskip+1pt\relax}%
    \def\thealgruledepth{\dimexpr .25\baselineskip+1pt\relax}%
}

\providecommand{\avg}[2]{\frac{1}{#2}\sum_{#1=1}^{#2}}

\usepackage{mathtools}
\DeclarePairedDelimiterX{\inp}[2]{\langle}{\rangle}{#1, #2}
\DeclarePairedDelimiterX{\cbr}[1]{\{}{\}}{#1} %
\DeclarePairedDelimiterX{\rbr}[1]{(}{)}{#1} %
\DeclarePairedDelimiterX{\sbr}[1]{[}{]}{#1} %

\usepackage{natbib}

\usepackage{thm-restate}

\usepackage{cleveref}

\theoremstyle{plain}
\newtheorem{theorem}{Theorem}[section]

\newtheorem{lemma}[theorem]{Lemma}

\theoremstyle{definition}
\newtheorem{definition}[theorem]{Definition}
\newtheorem{assumption}[theorem]{Assumption}
\newtheorem{notation}[theorem]{Notation}
\newtheorem{fact}[theorem]{Fact}
\theoremstyle{remark}
\newtheorem{remark}[theorem]{Remark}

\crefname{assumption}{Assumption}{Assumptions}
\crefname{fact}{Fact}{Facts}

\newcommand{\new}{\marginpar{NEW}}

\definecolor{darkgreen}{rgb}{0.0,0.5,0.0}
\definecolor{lightred}{rgb}{1, 0.8, 0.8}
\definecolor{lightgreen}{rgb}{0.8, 1, 0.8}

\newboolean{commentsOn}
\setboolean{commentsOn}{false}

\commenter{seb}{red}
\commenter{rick}{blue}
\commenter{anton}{orange}
\commenter{jerry}{green}

\usepackage[font=small]{caption}
\usepackage{cleveref}

\usepackage{import}

\crefname{equation}{}{}  %

\title{Accelerated Distributed Optimization with Compression and Error Feedback}

\author{
  Yuan Gao\thanks{CISPA Helmholtz Center for Information Security, Germany.} \thanks{Universit\"at des Saarlandes, Germany.} \\
    \texttt{yuan.gao@cispa.de}
    \and
    Anton Rodomanov\footnotemark[1]\\
    \texttt{anton.rodomanov@cispa.de}
    \and
    Jeremy Rack\footnotemark[1] \footnotemark[2] \\
    \texttt{michael.rack@cispa.de}
    \and
    Sebastian U. Stich\footnotemark[1] \\
    \texttt{stich@cispa.de}
}

\begin{document}

\maketitle

\begin{abstract}
	Modern machine learning tasks often involve massive datasets and models, necessitating distributed optimization algorithms with reduced communication overhead. Communication compression, where clients transmit compressed updates to a central server, has emerged as a key technique to mitigate communication bottlenecks. However, the theoretical understanding of stochastic distributed optimization with contractive compression remains limited, particularly in conjunction with Nesterov acceleration---a cornerstone for achieving faster convergence in optimization.
	In this paper, we propose a novel algorithm, \algname{ADEF} (\textbf{A}ccelerated \textbf{D}istributed \textbf{E}rror \textbf{F}eedback), which integrates Nesterov acceleration, contractive compression, error feedback, and gradient difference compression. We prove that \algname{ADEF} achieves the first accelerated convergence rate for stochastic distributed optimization with contractive compression in the general convex regime. Numerical experiments validate our theoretical findings and demonstrate the practical efficacy of \algname{ADEF} in reducing communication costs while maintaining fast convergence.
\end{abstract}

\section{Introduction}
\label{sec:introduction}
Gradient methods are the workhorses of the modern machine learning paradigm. They are the backbone for training modern machine learning models. With the rapidly growing sizes of the datasets and models, distributed training has become a necessity as accumulating the entire dataset on a single central machine is simply infeasible. Much of the recent breakthrough models are trained in such a distributed fashion, e.g.\ large language models \citep{shoeybi2019megatron}, generative models \citep{ramesh2021zero, ramesh2022hierarchical}, and others \citep{wang2020survey}, where data are distributed across different clients/workers and model updates are coordinated by a parameter server. Another instances of distributed optimization is the federated learning setting~\citep{FEDLEARN,kairouz2019advances}, where clients (e.g. edge devices or hospitals) jointly train a model without sharing their local data.

One of the key challenges in distributed optimization is the \emph{communication bottleneck}---the cost of transmitting large model updates between clients and the server \citep{seide20141,strom15_interspeech}.
To address this, \textbf{communication compression} has emerged as a practical solution, where updates are transmitted in a compressed format to reduce communication overhead.
Contractive compression (see \Cref{def:compression}), despite being potentially biased, has been shown to outperform unbiased compression schemes in practice \citep{lin2018deep,sun2019sparse,vogels2019PowerSGD} and offers favorable theoretical properties \citep{stich2018sparsified,albasyoni2020optimal,beznosikov2023biased}.
However, naive aggregation of compressed updates can lead to divergence \citep{beznosikov2023biased}. To mitigate this, \textbf{Error Feedback (EF)} \citep{seide20141} has been proposed as an effective remedy, compensating for compression errors \citep{stich2018sparsified,karimireddy2019error}. EF and its variants (e.g., \algname{PowerSGD} by \citet{vogels2019PowerSGD}) have been widely adopted in practice, integrated into popular deep learning frameworks such as PyTorch \citep{pazske2019pytorch}, and deployed to train state-of-the-art transformer models \citep{ramesh2021zero}.

Despite the empirical success of EF and contractive compression, their theoretical understanding remains incomplete. Early theoretical studies of EF \citep{stich2018sparsified,karimireddy2019error} primarily focus on single-client or homogeneous settings. Recent works have extended EF to heterogeneous datasets using bias correction mechanisms \citep{fatkhullin2023momentum,gao2024econtrol}, while the initial explorations often require additional unbiased compression or suffering from divergence in stochastic settings \citep{DIANA,richtarik2021ef21}.

On the other hand, most of these advances are made with unaccelerated methods. Nesterov acceleration~\citep{nesterov1983method} is one of the most important algorithmic tools for achieving faster convergence rates for deterministic and stochastic first order methods in the convex regime. 
However, its integration with EF and contractive compression in the general convex regime remains largely unexplored. Existing works on accelerated methods with compression primarily focus on deterministic or strongly convex settings \citep{li2020acceleration,qian2021error,bylinkin2024accelerated}, rely on unbiased compression \citep{li2021canita}, or effectively function as uncompressed methods \citep{he2023lower}.\footnote{We explain this in more details in \Cref{sec:history-acceleration} and \Cref{sec:repeated-communication}.} Consequently, the question of whether Nesterov acceleration can be achieved alongside contractive compression in the general convex regime remains open.

\subsection{Contribution} In this paper, we partially address the aforementioned gap by proposing \algname{ADEF} (\textbf{A}ccelerated \textbf{D}istributed \textbf{E}rror \textbf{F}eedback), which integrates Nesterov acceleration, contractive compression, and bias-corrected error feedback in the stochastic and general convex regime. Our contributions are as follows:

\begin{itemize}[topsep=0pt,itemsep=0pt,parsep=0pt]
\item \textbf{Algorithm:} We introduce \algname{ADEF}, a novel algorithm for accelerated distributed optimization that combines Nesterov acceleration with contractive compression, error feedback, and gradient difference compression. %
    
\item \textbf{General theoretical framework (acceleration with inexact updates):} We develop a general and versatile analysis framework for accelerated methods with inexact updates, which is of independent interest to the optimization and machine learning communities. This framework rigorously handles errors arising from compression, stochasticity, and other inexactness in distributed optimization. 
We demonstrate the flexibility of our framework by applying it to analyze other algorithms and compression techniques. %
This highlights the robustness of our approach and its potential to inspire future advances in distributed optimization.

\item \textbf{Breakthrough theoretical results:} Using our framework, we prove that \algname{ADEF} is the first distributed optimization algorithm to achieve accelerated convergence rates under contractive compression in the general convex regime. This result bridges a critical gap in the theoretical understanding of communication-efficient optimization and acceleration, overcoming longstanding challenges in the field. 

 \textbf{Numerical evaluation:} We conduct extensive numerical experiments to validate our theoretical findings. The results demonstrate that \algname{ADEF} significantly reduces communication costs while achieving state-of-the-art convergence performance.

\end{itemize}

By bridging the gap between communication-efficient distributed optimization and accelerated convergence, \algname{ADEF} opens new avenues for scalable and efficient training of large-scale machine learning models.

\section{Communication Compression and Accelerated Methods}

\begin{table*}[tb]

    \caption{Comparison of key existing works on stochastic distributed optimization with contractive compression in the general convex regime. We compare the rates in terms of the number of communication rounds from the client to the server, when the local oracle complexities are minized. For simplicity, we assume that $L=\ell=L_{\max}$.
    }
    \label{tab:comparison}
    \centering
    \resizebox{\linewidth}{!}{
    \begin{tabular}{ccccc}
        \toprule
        \textbf{Algorithm} & $\textbf{Acc}^{\rm (a)}$  & $\textbf{No BGS}^{\rm (b)}$ & \# Comm. Rounds & \# Local Oracle Calls \\
        \midrule
        \makecellnew{Lower Bound\\\citep{he2023lower}} & - & - & $\Omega\left(\frac{R_0^2\sigma^2}{n\varepsilon^2} + \frac{\sqrt{LR_0^2}}{\delta\sqrt{\varepsilon}}\right)$& $\Omega\left(\frac{R_0^2\sigma^2}{n\varepsilon^2} + \frac{\sqrt{LR_0^2}}{\sqrt{\varepsilon}}\right)$\\
        \cmidrule{1-5}
        \makecellnew{\algname{EF}\\\citep{karimireddy2019error}} & No & No & $\cO\left(\frac{R_0^2\sigma^2}{n\varepsilon^2}+\frac{\sqrt{L}R_0^2(\frac{\zeta}{\sqrt{\delta}}+\sigma)}{\sqrt{\delta}\varepsilon^{\nicefrac{3}{2}}}+\frac{LR_0^2}{\delta\varepsilon}\right)$ & $\cO\left(\frac{R_0^2\sigma^2}{n\varepsilon^2}+\frac{\sqrt{L}R_0^2(\frac{\zeta}{\sqrt{\delta}}+\sigma)}{\sqrt{\delta}\varepsilon^{\nicefrac{3}{2}}}+\frac{LR_0^2}{\delta\varepsilon}\right)$\\
        \cmidrule{1-5}
        \makecellnew{\algname{EControl}\\\citep{gao2024econtrol}} & No & Yes & $\cO\left(\frac{R_0^2\sigma^2}{n\varepsilon^2}+\frac{\sqrt{L}R_0^2\sigma}{\delta^2\varepsilon^{\nicefrac{3}{2}}}+\frac{LR_0^2}{\delta\varepsilon}\right)$& $\cO\left(\frac{R_0^2\sigma^2}{n\varepsilon^2}+\frac{\sqrt{L}R_0^2\sigma}{\delta^2\varepsilon^{\nicefrac{3}{2}}}+\frac{LR_0^2}{\delta\varepsilon}\right)$\\
        \cmidrule{1-5}
        \makecellnew{$\text{\algname{NEOLITHIC}}^{\rm (c)}$\\\citep{he2023lower}} & Yes & No & $\wtilde\cO\left(\frac{R_0^2\sigma^2}{\delta n\varepsilon^2}+\frac{\sqrt{LR_0^2}}{\delta\sqrt{\varepsilon}}\right)^{\rm (d)}$& $\cO\left(\frac{R_0^2\sigma^2}{ n\varepsilon^2}+\frac{\sqrt{LR_0^2}}{\sqrt{\varepsilon}}\right)$\\
        \cmidrule{1-5}
        \makecellnew{\algname{ADEF}\\new} & Yes & Yes & $\cO\left(\frac{R_0^2\sigma^2}{n\varepsilon^2}+ \frac{\sqrt{L}R_0^2\sigma}{\delta^2\varepsilon^{\nicefrac{3}{2}}}+\frac{\sqrt{LR_0^2}}{\delta^2\sqrt{\varepsilon}} \right)$ & $\cO\left(\frac{R_0^2\sigma^2}{n\varepsilon^2}+ \frac{\sqrt{L}R_0^2\sigma}{\delta^2\varepsilon^{\nicefrac{3}{2}}}+\frac{\sqrt{LR_0^2}}{\delta^2\sqrt{\varepsilon}} \right)$\\
        \bottomrule
    \end{tabular}
    }
    \vskip0.5em
    \begin{tablenotes}
        {\scriptsize
        \item (a) Acc stands for accelerated methods. 
        \item (b) BGS stands for Bounded Gradient Similarity assumption: $\avg{i}{n}\norm{\nabla f_i(\xx)-\nabla f(\xx)}^2\leq \zeta^2,\forall \xx\in\R^d$. See also \Cref{assumption:gradient-similarity}.
        \item (c) \citet{he2023lower} presented their rates when batch size is $\max\{\frac{4}{\delta}\ln(\frac{4}{\delta}), \frac{1}{\delta}\ln(24\kappa + \frac{25n^2\kappa^3\zeta^4}{\sigma^4}+5n\kappa^{\nicefrac{3}{2}})\}$ where $\kappa$ is the total number of updates at the server. We present the rate, as is the same with other methods, when the local oracle complexity is minimized, i.e. when the batch size is one. We also note that, in the worst case, \algname{NEOLITHIC} transmit more bits than simply sending the uncompressed vector between each update at the server and is therefore an uncompressed method.
        \item (d) $\wtilde\cO$ hides logarithmic factors in $\nicefrac{1}{\varepsilon}$.
        }
    \end{tablenotes}
\end{table*}

In this paper, we study distributed optimization algorithms equipped with communication compression. In particular, we consider the class of contractive compression operators:
\begin{definition}
    \label{def:compression}
    We say that a (possibly randomized) mapping  $\cC \colon \R^d \to \R^d$ is a contractive compression operator if for some constant $0 < \delta \le 1$ it holds 
    \begin{equation}
        \label{eq:compressor}
        \Eb{\norm{\cC(\xx) -\xx}^2} \leq (1-\delta)\norm{\xx}^2\, \quad \forall \xx \in \R^d. 
    \end{equation}
\end{definition}
The classic example satisfying \Cref{def:compression} is Top-$K$ \citep{stich2018sparsified}. Top-$K$ keeps the K largest entries in terms of absolute value and zeros out the rest and it is a biased compression. The class of contractive compressors also includes sparsification~\citep{alistarh2018convergence,stich2018sparsified} and quantization operators~\citep{wen2017terngrad,QSGD,bernstein2018signsgd,horvath2019natural}, and many others \citep{zheng2019communication,beznosikov2020biased, vogels2019PowerSGD,safaryanFedNL2022,islamov2023distributed}. We follow the assumption made in most existing works on the theory of communication compressed optimization that the server-to-client communication is much more efficient (as compared to client to server communication) and the broadcast cost is negligible~\citep{DIANA,kairouz2019advances,fatkhullin2023momentum,gao2024econtrol}. Therefore, we only consider client-to-server communication compression in this paper.

In \Cref{tab:comparison} we compare our method to some key existing works most relevant to our setting, stochastic distributed optimization with contractive compression in the general convex regime. We measure the complexity in terms of the number of communication rounds from the client to the server when the local oracle complexity is minimized, in other word, when the batch size is 1.

\subsection{Brief History on the Theory of \algname{EF}}
The \algname{EF} mechanism, proposed by \citet{seide20141}, was first analyzed theoretically in \citet{stich2018sparsified,alistarh2018convergence,karimireddy2019error}, but only for the single client setting. Extensions to the distributed setting were first made under data similarity assumptions, either implicitly in the form of bounded gradient assumption~~\citep{cordonnier2018convex,alistarh2018convergence}, or explicitly in the form of gradient similarity assumptions~\citep{stich2020error}, both are very limiting factors. Further extensions to fully decentralized settings were also considered in~\citep{koloskova2019Decentralized,koloskova2020Decentralized} under the bounded gradient assumption. The theory of distributed \algname{EF} were further refined in~\citep{beznosikov2020biased,stich2020communication}. A key point in the analysis of distributed \algname{EF} is to obtain an convergence rate in the number of comunication rounds, where the leading term (term that involves the variance of the stochastic oracle, see \Cref{assumption:bounded_variance}) are unaffected by the compression quality $\delta$ and enjoys the linear speedup in the number of clients.

\citet{DIANA} proposed the \algname{DIANA} algorithm, which incorporates an additional \emph{unbiased} compressor into \algname{EF} for bias correction, alleviating the need for the data similarity assumptions. It inspired a number of follow up works~\citep{gorbunov2020linearly,stich2020communication,qian2021errorcompensatedSGD}, and eventually led to the \algname{EF21} algorithm, which is the first one that fully supports contractive compression in the full gradient regime~\citep{richtarik2021ef21}. However, the bias correction meschanism of \algname{EF21} does not work with stochastic gradients and leads to unconvergence up to the variance of the stochastic oracle. The challenge of bias correction with stochastic gradients was finally addressed in~\citet{fatkhullin2023momentum} with the application of momentum (which however, is not known to work in the general convex regime), and in~\citet{gao2024econtrol} using an error-controlled \algname{EF} combined with bias correction which covers the strongly convex, general convex and non-convex cases.

\subsection{Accelerated Methods with Compression}
\label{sec:history-acceleration}
Compared to the rich literature on the unaccelerated \algname{EF} and its variants, theory for accelerated methods~\citep{nesterov1983method} with communication compression is still lacking, with only a handful available~\citep{murata2019accelerated,li2020acceleration,qian2021error,li2021canita,he2023lower,bylinkin2024accelerated}. This is particularly the case with contractive compression, where some work~\citep{he2023lower} see the potential biasedness in the updates a forbiding factor, due to the negative results~\citep{devolder2014first} on acceleration with inexact oracles. %

\textbf{Unbiased Compression:} \algname{ADIANA} \citep{li2020acceleration} and \algname{CANITA} \citep{li2021canita} are both restricted to unbiased compressions. %
The \algname{S-NAG-EF} proposed in~\citep{murata2019accelerated} even requires a more specific class of randomized unbiased compressors. We also note that a very recent work by~\citet{gupta2024nesterov} essentially rediscovered the theory for compressed \algname{GD} with unbiased compressors and its accelerated variant~\citep{li2020acceleration}, under the notion of noisy gradient with relative noise. 
The unbiasedness assumption is a very strong assumption, and it is known in practice that contractive compresors are superior~\citep{vogels2019PowerSGD,albasyoni2020optimal,beznosikov2023biased}, rendering these theories unsatisfactory.

\textbf{Lack of Theory for General Convex Case:} \algname{ADIANA}, \algname{ECLK} \citep{qian2021error} and \algname{EF-OLGA} \citep{bylinkin2024accelerated} only work for the strongly-convex regime. It is not known if these methods can be extended to the general convex case.

\textbf{Only Supports Full Gradient:} As we have seen before, mixing stochastic gradients with compression can already be problematic in the unaccelerated case. \algname{ADIANA}, \algname{ECLK}, \algname{CANITA} and \algname{EF-OLGA} all only work with full gradients. In particular, \algname{ECLK} and \algname{EF-OLGA} combines \algname{EF} with \algname{Katyusha}-style accceleration~\citep{allen2018katyusha}, which utilizes the special finite-sum structure of $f$ and requires epoch-wise \emph{uncompressed} communications of \emph{full} gradients, making it unlikely to be extendable to the stochastic setting. \algname{ADIANA} and \algname{CANITA} does not utilize the classic \algname{EF} mechanism, and evidences suggest that their gradient difference compression mechanism might lead to unconvergence up to the variance of the stochastic oracle when combined with contractive compression~\citep{fatkhullin2023momentum}.

\textbf{Uncompressed Method in Disguise:} \algname{NEOLITHIC} \citep{huang2022lowerbounds,he2023lower} takes a rather different approach than the above. %
It does not utilize any additional algorithmic tools on top of a naively compressed Nesterov accelerated method. It uses $\wtilde\Omega(\frac{1}{\delta})$ rounds\footnote{More precisely: $R=\max\{\frac{4}{\delta}\ln(\frac{4}{\delta}), \frac{1}{\delta}\ln(24\kappa + \frac{25n^2\kappa^3\zeta^4}{\sigma^4}+5n\kappa^{\nicefrac{3}{2}})\}$. Here $\zeta^2$ denotes the gradient similarity and $\kappa$ denotes the total number of updates} of communications between each update to account for compression errors. However, this means that in the worst case the algorithm will transmit \emph{more} bits than simply sending the uncompressed vector between each server updates, making it an \emph{uncompressed} method in disguise. 
See \Cref{sec:repeated-communication} for a more in-depth discussion on these issues.

\section{Problem Formulation and Assumptions}
\label{sec:problem_formulation}
In this paper we consider the following distributed stochastic optimization problem:
\begin{equation}
    \label{eq:problem}
    \xx^\star \eqdef \argmin_{\xx\in\R^d} \Bigl[ f(x) \eqdef \avg{i}{n}f_i(\xx) \Bigr],
\end{equation}
where $\xx$ are the parameters of a model that we train. The objective function $f:\R^d\to\R$ is an average of $n$ functions $f_i:\R^d\to\R,i\in[n]$. Each function $f_i$ is a local loss associated with a local dataset $\cD_i$ which can only be accessed by client $i$.

We analyze the setting where there is a central server coordinating the training process. The server receices and aggregates the messages from the clients via a compressed communication channel, where the messages are compressed by some $\delta$-contractive compression (see \Cref{def:compression}), performs an update on the model parameters, and broadcasts the updated model back to the clients. We measure the performance of the algorithm in terms of the number of communication rounds from the clients to the server, when the local oracle complexity of the clients are minimized.

We first introduce the following notation which we use throughout the paper.

\begin{notation}
    \label{notation:divergence}
    We denote the Bregman divergence of $f$ between $\xx$ and $\yy$ as $\beta_f(\xx,\yy)$, i.e.,
    \begin{equation}
        \label{eq:divergence}
        \beta_f(\xx,\yy) \eqdef f(\yy)-f(\xx)-\inp{\nabla f(\xx)}{\yy-\xx}\,.
    \end{equation}
\end{notation}

Next we list the assumptions that we make in this paper. Notably, we do not assume that the local functions $f_i$'s have gradient dissimilarity, which is a common limiting assumption in many of the existing works, in both the unaccelerated and accelerated cases, e.g. \citep{lian2017can,huang2022lowerbounds,li2023analysis,he2023lower}. A crucial algorithmic component of our method that enables the acceleration also achieves this bias correction pheonomenon naturally.

\begin{assumption}
    \label{assumption:convexity}
    We assume that the objective function $f$ is convex, closed and proper.
\end{assumption}

\begin{assumption}
    \label{assumption:smoothness}
    We assume that the objective function $f$ has $L$-Lipschitz gradients, i.e. for all $\xx,\yy\in\R^d$, it holds
    \begin{equation}
        \label{eq:smoothness}
        \norm{\nabla f(\xx) - \nabla f(\yy)} \leq L\norm{\xx-\yy}.
    \end{equation}
\end{assumption}

Many existing works further assumes that each local function $f_i$ is $L_{\max}$-smooth and convex~\citep{li2021canita,he2023lower}. We instead make the following weaker assumption:
\begin{assumption}
    \label{assumption:f-i-smoothness}
    We assume that there exists some $\ell>0$ such that for all $\xx,\yy\in\R^d$, it holds
    \begin{equation}
        \label{eq:f-i-smoothness}
        \avg{i}{n}\norm{\nabla f_i(\xx) - \nabla f_i(\yy)}^2 \leq 2\ell \beta_f(\xx,\yy).
    \end{equation}
\end{assumption}
It is easy to show that \Cref{assumption:f-i-smoothness} follows from $L_{\max}$-smoothness of each local function $f_i$, in particular:
\begin{fact}
    \label{fact:f-i-smoothness}
    Given the assumption that each $f_i$ is $L_{\max}$-smooth and convex, then \Cref{assumption:f-i-smoothness} holds with $\ell\leq L_{\max}$.
\end{fact}
Another class of common assumptions used in the literature are variants of the Hessian similarity assumption~\citep{khaled2022faster,jiang2024federated,rodomanov2024universality,bylinkin2024accelerated}, which also implies \Cref{assumption:f-i-smoothness}:
\begin{fact}
    \label{fact:hessian_similarity-f-i-smoothness}
    Given \Cref{assumption:convexity,assumption:smoothness}, if the local functions $f_i$'s satisfy the $\lambda$-Hessian similarity assumption in the following sense:
    \begin{equation}
        \label{eq:hessian_similarity-f-i-smoothness}
        \avg{i}{n}\norm{\nabla \hat f_i(\xx)-\nabla \hat f_i(\yy)}^2 \leq 2\lambda \beta_f(\xx,\yy)\,,
    \end{equation}
    where $\hat f_i\eqdef f-f_i$, then $\ell\leq L+\lambda$.
\end{fact}

Next we make the following standard assumption on the stochastic gradient oracles.
\begin{assumption}
    \label{assumption:bounded_variance}
    We assume that each client $i$ has access to an unbiased stochastic gradient oracle $\gg_i(\xx,\xi^i):\R^d\to\R^d$ for the local function $f_i$, such that for all $\xx\in\R^d$, it holds
    \begin{equation}
        \label{eq:bounded_variance}
        \begin{aligned}
            &\EEb{\xi^i}{\gg_i(\xx,\xi^i)} = \nabla f_i(\xx),\\
            &\EEb{\xi^i}{\norm{\gg_i(\xx,\xi^i)-\nabla f_i(\xx)}^2} \leq \sigma^2.
        \end{aligned}
    \end{equation}
\end{assumption}
Mini-batches are also allowed, which simply divides the variance by the batch size. Our algorithm does not impose restrictions on the minimal batch size needed during training. We compare the total number of communiation rounds from the client to the server for the algorithms, when the loca oracle complexity (i.e. the numbr of oracle calls of the client) is minimized, i.e. when the batch size is 1. This is also the standard metric in most existing works on stochastic distributed optimization with compression~\citep{karimireddy2019error,stich2020communication,fatkhullin2023momentum,gao2024econtrol}.

\section{Accelerated Method with Inexact Update}
\label{sec:acc-inexact}

\begin{algorithm}[tb]
    \caption{Accelerated Method with Inexact Update}
    \label{alg:acc-inexact}
    \begin{algorithmic}[1]
        \State \textbf{Input:} $\xx_0, \vv_0, A_0, (a_t)_{t=1}^\infty$
        \For{$t = 0,1,\dots$}
        \State $A_{t+1} = A_t+a_{t+1}$
        \State $\yy_t = \frac{A_t}{A_{t+1}}\xx_t+ \frac{a_{t+1}}{A_{t+1}}\vv_t$
        \State  Compute $\hat \gg_t \approx \gg_t$ %
        \State $\vv_{t+1} = \vv_t -a_{t+1}\hat\gg_t$
        \State $\xx_{t+1} = \frac{A_t}{A_{t+1}}\xx_t + \frac{a_{t+1}}{A_{t+1}}\vv_{t+1}$
        \hfill
        \EndFor
    \end{algorithmic}
  \end{algorithm}

In this section, we present a general framework for studying accelerated methods with inexact updates. Here we do not make any assumptions on the inexactness and obtain a general framework for the analysis. Our framework is flexible and can be applied to different setups and algorithms, unlike some existing works such as \citet{devolder2014first} which assumes specific structures on the inexactness. We consider \Cref{alg:acc-inexact}. The algorithm is a simple extension to the standard accelerated method, where in Line 5 we use some inexact $\hat\gg_t$ that is ``approximately'' the stochastic gradient $\gg_t$ of $f$ at $\yy_t$. We assume the following for the stochastic gradient $\gg_t$:

\begin{assumption}
    \label{assumption:bounded_variance_general}
    For all $t\geq 0, \gg_t\eqdef \gg(\yy_t,\xi_t)$ is a stochastic gradient of $f$ at $\yy_t$,
    where $\xi_t$ are independent copies of the oracle's randomness~$\xi$.
    We assume that for $\xx\in\R^d, \gg(\xx,\xi)$ is an unbiased stochastic gradient with bound variance:
    \begin{equation}
        \label{eq:bounded_variance_general}
        \begin{aligned}
            &\Eb{\gg(\xx,\xi)} = \nabla f(\xx),\\
            &\EEb{\xi}{\norm{\gg(\xx, \xi)-\nabla f(\xx)}^2} \leq \sigma_{\gg}^2.
        \end{aligned}
    \end{equation}
    Furthermore, we assume that $\hat\gg_t$ is a deterministic function of $\xi_{[t]}\eqdef \{\xi_0,\ldots,\xi_{t}\}$.
\end{assumption}

The value of $\sigma^2_{\gg}$ can be adapted to different setups, and in particular, for distributed optimization problems, we think of $\gg_t$ as the average of the local stochastic gradients (see \Cref{assumption:bounded_variance}) and therefore we have $\sigma^2_{\gg}=\frac{\sigma^2}{n}$. For now, we do not assume any specific structure on the inexact gradient $\hat\gg_t$. Such an inexact gradient could be the result of absolute compression (see \Cref{sec:absolute-compression}); or it could be built algorithmically via some compressed communication channel, which is the focus of this paper. 

The main result in this section is a descent theorem for the final iterate of~\Cref{alg:acc-inexact}, up to some accumulative error that depends on the inexactness of the gradient $\hat\gg_t$. In particular, we define the following accumulative error at iteration $t$:
\begin{equation}
    \label{eq:error-definition}
    \ee_t \eqdef \sum_{j=0}^{t-1}a_{j+1}(\hat \gg_j-\gg_j).
\end{equation}
We also introduce some additional notations:
\begin{equation}
    \label{eq:quantity-definition}
    \begin{aligned}
        F_t &\eqdef \E{f(\xx_t)-f(\xx^\star)},\\
        R_t^2 &\eqdef \EbNSq{\vv_t-\xx^\star},\\
        E_t &\eqdef \EbNSq{\ee_t}\,. %
    \end{aligned}
\end{equation}
Now, the main result of this section is the following:
\begin{restatable}{theorem}{mainDescent}
    \label{thm:main-descent}
    Given
    \Cref{assumption:convexity,assumption:smoothness,assumption:bounded_variance_general},
    for all $T\geq 1$, for $(\xx_t,\yy_t,\vv_t)_{t=0}^\infty$ generated by \Cref{alg:acc-inexact}, it holds that:
    \begin{equation}
        \label{eq:main-descent}
        \begin{aligned}
            A_TF_T &\leq A_0F_0+\frac{R_0^2}{2} + 2\sigma^2_g\sum_{t=0}^{T-1}a_{t+1}^2 + \sum_{t=0}^{T-1}w_{t+1}E_{t+1}\\
            &\quad  -\sum_{t=0}^{T-1}\left(\frac{1}{6}-\frac{2La_{t+1}^2}{A_{t+1}}\right)a_{t+1}^2\EbNSq{\gg_t}\\
            &\quad -\sum_{t=0}^{T-1}\left( \Eb{\frac{a_{t+1}}{2}\beta_f(\yy_t,\xx^\star) + A_t\beta_f(\yy_t,\xx_t) }\right) \\
            &\quad - \sum_{t=0}^{T-1}\Eb{A_{t+1}\beta_f(\xx_{t+1},\yy_t)}
        \end{aligned}
    \end{equation}
    where we write $w_t\eqdef \min\{2,a_tL\}+\frac{4La_t^2}{A_t}+\frac{4La_{t+1}^2}{A_{t+1}}$.
\end{restatable}

The missing proof is in \Cref{sec:acc-inexact-proof}. We point out that while \Cref{thm:main-descent}'s proof employs virtual iteration techniques used in the analysis of \algname{EF} and its variants (e.g. \citep{karimireddy2019error,stich2020communication,gao2024econtrol}), the proof is more involved as it also carefully integrates the acceleration scheme. More importantly, \Cref{thm:main-descent} takes a different appoach from the one-step descent analysis template which most existing works applied~\citep{karimireddy2019error,stich2020communication,richtarik2021ef21,fatkhullin2023momentum,gao2024econtrol}, in that it is an overall descent lemma for the final iterate. Instead of constructing a single Lyapunov function that attempts to captures all the indivual one-step descent (e.g. $R_t^2$ and $E_t$ etc) that in many cases involves quite some speculations for the appropriate choices of parameters, we only need to upper bound the sum of the accumulative errors now.

As a technicality in the analysis, we weight each error term $E_{t+1}$ with $\min\{2M,a_{t+1}L\}$, instead of simply using either one throughout. This comes from the fact that the stepsize $a_{t+1}$, albeit typically growing linearly, might attain lower value than a fixed constant throughout the algorithm due to the stepsize tunning for the stochastic terms.

As discussed before, \Cref{thm:main-descent} covers a wide range of algorithms that takes the form of \Cref{alg:acc-inexact} and can be applied to study the convergence of these algorithms as long as the sum of the accumulative errors can be controlled. To demonstrate its flexibility, we apply \Cref{thm:main-descent} to analyze various algorithms besides our main algorithm in this paper:
\begin{itemize}
    \item in \Cref{sec:acceleration-vanilla-ef}, we analyze \Cref{alg:dist-acc-ef-vanilla} where $\hat\gg_t$ is built via the vanilla error feedback, and demonstrate that it does not achieve the accelerated rate;
    \vspace{-0.5em}
    \item in \Cref{sec:repeated-communication}, we analyze \Cref{alg:neolithic}, which is essentially the \algname{NEOLITHIC} algorithm in~\citet{he2023lower}, where $\hat\gg_t$ is built via repeated communication rounds.
    \vspace{-0.5em}
    \item in \Cref{sec:absolute-compression}, we briefly study a different class of compressors, absolute compressor $\cC_{\Delta}$ such that $\Eb{\norm{\cC_{\Delta}(\xx)-\xx}^2}\leq \Delta^2$, and analyze the convergence of \Cref{alg:dist-acc-absolute} with such compressors.
\end{itemize}
\vspace{-0.5em}
Most importantly, in the next section, we present our main algorithm for distributed optimization with communication compression, \algname{ADEF}, and analyze its convergence using \Cref{thm:main-descent}.

\section{Acceleration with Error Feedback}
\label{sec:acc-ef}
\begin{algorithm}[tb]
    \caption{\algname{ADEF} Accelerated Distributed Error Feedback}
    \label{alg:dist-acc-ef}
    \begin{algorithmic}[1]
        \State \textbf{Input:} $\xx_0, \vv_0, A_0, (a_t)_{t = 1}^\infty$ and $\cC$.
        \State $\ee_0^i=\0$, \ $\tilde{\gg}_{-1}^i = \gg_i(\yy_0,\xi_0^i)$, \ $\tilde{\gg}_{-1} = \avg{i}{n}\wtilde\gg_{-1}^i$
        \For{$t = 0,1,\dots$}
        \State {\bf server:}
        \State $A_{t+1} = A_t+a_{t+1}$
        \State send to each client $\yy_t = \frac{A_t}{A_{t+1}}\xx_t+ \frac{a_{t+1}}{A_{t+1}}\vv_t$
        \State {\bf each client $i$:}
        \State $\gg_t^i = \gg_i(\yy_t,\xi_t^i)$
        \State \colorbox{lightred}{$\tilde{\delta}_t^i=\gg_t^i-\tilde{\gg}_{t-1}^i$, \ $\tilde{\Delta}_t^i = \cC( \tilde{\delta}_t^i)$}
        \State \colorbox{lightred}{$\tilde{\gg}_t^i=\tilde{\gg}_{t-1}^i+\tilde{\Delta}_{t}^i$}
        \State \colorbox{lightgreen}{$\delta_{t}^i= \gg_t^i-\tilde{\gg}_{t}^i - \frac{1}{a_{t+1}} \ee_t^i$, \ $\Delta_{t}^i = \cC(\delta_{t}^i)$}
        \State \colorbox{lightgreen}{$\ee_{t+1}^i = a_{t+1}( \Delta_{t}^i - \delta_{t}^i)$}
        \State send to server $\tilde{\Delta}_t^i$, $\Delta_{t}^i$
            \State {\bf server:}
        \State \colorbox{lightred}{$\tilde{\gg}_t = \tilde{\gg}_{t-1} + \frac{1}{n} \sum_{i=1}^{n}\tilde{\Delta}_{t}^i$}
        \State \colorbox{lightgreen}{$\hat \gg_t = \tilde{\gg}_t + \frac{1}{n} \sum_{i=1}^{n}\Delta_{t}^i$}
        \State $\vv_{t+1}= \vv_t-a_{t+1}\hat \gg_t$
        \State $\xx_{t+1} = \frac{A_t}{A_{t+1}}\xx_t + \frac{a_{t+1}}{A_{t+1}}\vv_{t+1}$
        \EndFor
    \end{algorithmic}	
\end{algorithm}

In this section we present our main algoirthm, \algname{ADEF}, summarized in \Cref{alg:dist-acc-ef}. All missing proofs can be found in \Cref{sec:missing-proofs}. At its core, the algorithm is still the accelerated method with inexact updates in the form of \Cref{alg:acc-inexact}. However, we introduce an enhanced error feedback mechanism to build the inexact gradient $\hat \gg_t$, and we highlight the procedures in which we build the inexact gradient with red and green color boxes in \Cref{alg:dist-acc-ef}. The algorithm employs a delicate combination of error feedback and gradient difference compression. The green lines, Line 11, 12, 16, is the backbone error feedback procedure:
\begin{equation}
    \label{eq:ef-procedure}
    \begin{aligned}
        \delta_{t}^i&= \gg_t^i-\tilde{\gg}_{t}^i - \frac{1}{a_{t+1}} \ee_t^i, \quad \Delta_{t}^i = \cC(\delta_{t}^i),\\
        \ee_{t+1}^i &= a_{t+1}( \Delta_{t}^i - \delta_{t}^i),\quad \hat{\gg}_t = \tilde{\gg}_t + \frac{1}{n} \sum_{i=1}^{n}\Delta_{t}^i.
    \end{aligned}
\end{equation}
For any sequences $\wtilde\gg_t^i$, the error feedback mechanism summarized above records the errors that occured due to the compression. The following lemma, which is the key property of \algname{EF} and most of its variants~\citep{karimireddy2019error,stich2020communication,gao2024econtrol}, states that the accumulative error defined in \Cref{eq:error-definition} is precisely the average of the local errors tracked by \Cref{eq:ef-procedure}, irrespective of the choice of $\wtilde\gg_t^i$'s. Consequently, we can bound the sum of the accumulative errors $\sum_{t=0}^{T-1}w_{t+1}E_{t+1}$ with the following quantity that quantifies the distance between the local stochastic gradients and the control variates:
\begin{equation}
    \label{eq:quantity-definition-h}
    H_t \eqdef \avg{i}{n}\EbNSq{\gg_t ^i - \wtilde\gg_t ^i}.
\end{equation}

Now we have:
\begin{restatable}{lemma}{eDescentSumInH}
    \label{lem:e-descent-sum-in-h}
    For an algorithm that follows \Cref{eq:ef-procedure}, the local error terms $(\ee_T^i)_{t=0}^\infty$ satisfies the following:
    \begin{equation}
        \label{eq:error-identity}
        \avg{i}{n}\ee_t^i = \sum_{j=0}^{t-1}a_{j+1}(\hat \gg_j-\gg_j)\eqcolon \ee_t,\,
    \end{equation}
    where $\gg_t$'s are the average local stochastic gradients $\gg_t = \avg{i}{n}\gg_t^i$.
    
    Further, if for all $t\leq T-1, (1-\frac{\delta}{2})w_{t+1}\leq (1-\frac{\delta}{4})w_t$, then it holds that:
    \begin{equation}
        \label{eq:e-descent-sum-in-h}
            \sum_{t=0}^{T-1}w_{t+1}E_{t+1} \leq \frac{8}{\gamma}\sum_{t=0}^{T-1}w_{t+1}'H_t,\,
    \end{equation}
    where $w_t'\eqdef w_ta_t^2$
\end{restatable}
Now, for the vanilla error feedback method, the control variate $\wtilde\gg_t^i$'s are simply $\0$. This is however undesirable, as upper bounding $w_{t+1}'H_t$ will involve upper bounding $\frac{w_{t+1}'}{n}\sum_{i=1}^{n}\EbNSq{\nabla f_i(\yy_t)}$. This in turns will require the data similarity assumption, and more importantly, an upper bound in the form of $2w_{t+1}'\ell \beta_f(\yy_t,\xx^\star)$. This cannot be canceled out unless we force the stepsizes to be very small which effectively eliminate all the acceleration. See \Cref{sec:acceleration-vanilla-ef} for more details. To address these issues with the vanilla \algname{EF} mechanism, we build the control variate $\wtilde\gg_t^i$'s to approximate $\gg_t^i$'s as well and therefore reduce the overall errors from the compression. The construction of $\wtilde\gg_t^i$'s are highlighted in red in \Cref{alg:dist-acc-ef} and we summarize it below:
\begin{equation}
    \label{eq:control-variate-procedure}
    \begin{aligned}
        \gg_t^i &= \gg_i(\yy_t,\xi_t^i), \quad \tilde{\delta}_t^i= \gg_t^i-\tilde{\gg}_{t-1}^i, \\
        \tilde{\Delta}_t^i &= \cC( \tilde{\delta}_t^i),\quad \tilde{\gg}_t^i =\tilde{\gg}_{t-1}^i+\tilde{\Delta}_{t}^i.
    \end{aligned}
\end{equation}
The procedure described in \Cref{eq:control-variate-procedure} is sometimes known as a gradient difference compression mechanism, which are mostly applied to address the data heterogeneity issue for unaccelerated methods~\citep{DIANA,fatkhullin2021ef21,gao2024econtrol}. Here, it reduces the error of compression and addresses the daaa heterogeneity issue all at once. We have the following lemma regarding the weighted sum of $H_t$:

\begin{restatable}{lemma}{hDescentSum}
    \label{lem:h-descent-sum}
    Given \Cref{assumption:f-i-smoothness,assumption:bounded_variance}, for an algorithm that follows \Cref{eq:control-variate-procedure}, if for all $t\leq T-1$, we have $w_{t+1}'(1-\frac{\delta}{2}) \leq w_t'(1-\frac{\delta}{4})$, then it holds that:
    \begin{equation}
        \label{eq:h-descent-sum}
        \begin{aligned}
            \sum_{t=0}^{T-1}w_{t+1}'H_t & \leq  \frac{16}{\delta^2}\sum_{t=1}^{T-1}w_{t+1}'\avg{i}{n}\EbNSq{\nabla f_i(\yy_t)-\nabla f_i(\yy_{t-1})} + \frac{16\sigma^2}{\delta^2}\sum_{t=1}^{T-1}w_{t+1}'.
        \end{aligned}
    \end{equation}
\end{restatable}

Finally, we can combine \Cref{thm:main-descent} with \Cref{lem:e-descent-sum-in-h} and \Cref{lem:h-descent-sum} to give the main convergence result for \Cref{alg:dist-acc-ef}. As we pointed out before, for the distributed optimization setting that we consider in this section, $\gg_t$ is simply the average of the local gradient and therefore $\sigma^2_{\gg}=\frac{\sigma^2}{n}$. 

\begin{restatable}{theorem}{convergenceRate}
    \label{thm:convergence-rate}
    Given \Cref{assumption:convexity,assumption:smoothness,assumption:f-i-smoothness,assumption:bounded_variance}, and let $a_t\eqdef \nicefrac{(t+\frac{32}{\delta})}{M}$ and $A_0\eqdef \nicefrac{32^2}{2\delta^2 M}$, it suffices to have:
    \begin{equation}
        \label{eq:convergence-rate}
        T = \cO\left(\frac{R_0^2\sigma^2}{n\varepsilon^2}+ \frac{\sqrt{L}R_0^2\sigma}{\delta^2\varepsilon^{\nicefrac{3}{2}}}+\frac{\sqrt{\ell R_0^2}}{\delta^2\sqrt{\varepsilon}} \right),\,
    \end{equation}
    number of iterations (communication rounds) of \Cref{alg:dist-acc-ef} to get $F_T\leq \varepsilon$, where we can set $M=\max\left\{\frac{2^{13}\ell}{\delta^4}, \left(\frac{4T(T+\frac{32}{\delta})^2\sigma^2}{R_0^2n}\right)^{\frac{1}{2}}, 8\left(\frac{LT(T+\frac{32}{\delta})^3\sigma^2}{\delta^4R_0^2}\right)^{\frac{1}{3}}\right\}$.
\end{restatable}
\begin{remark}
    We note that for the proof of \Cref{thm:convergence-rate}, we mainly focus on the asymptotic behavior and does not optimize the choices of constants. Most likely one can pick much smaller constants for $M$ in practice, see also \Cref{sec:experiments}.
\end{remark}
\begin{remark}
    \Cref{thm:convergence-rate} gives the first accelerated convergence rate for distributed stochastic optimization with error feedback and contractive compression. We first point out that the stochastic terms $\cO(\frac{R_0^2\sigma^2}{n\varepsilon^2}+ \frac{\sqrt{L}R_0^2\sigma}{\delta^2\varepsilon^{\nicefrac{3}{2}}})$ matches that of \algname{EControl}~\citep{gao2024econtrol}, which is the best known for the distributed setting under arbitrary data heterogeneity. In particular, the leading term $\cO(\frac{R_0^2\sigma^2}{n\varepsilon^2})$ is optimal, $\delta$-free, and enjoys the linear speedup in the number of clients $n$. 
    
    Importantly, in the deterministic regime, i.e. when $\sigma^2=0$, we see that \algname{ADEF} achieves the accelerated $\cO(\frac{1}{\sqrt{\varepsilon}})$ convergence rate. There is however a $\frac{1}{\delta^2}$ factor, instead of a $\frac{1}{\delta}$ factor. Nonetheless, this is faster than the unaccelerated $\cO(\frac{LR_0^2}{\delta \varepsilon})$ when $\delta\geq \sqrt{\frac{\varepsilon}{LR_0^2}}$ (for simplicity of the illustration, we assume that $L=\ell=L_{\max}$ in this comparison). We also point out that this higher dependence on the $\delta$ factor in the deterministic term seems to be a common issue in the analysis of bias-corrected \algname{EF} when using two compressions. In particular, such a pheonomenon was also observed in the analysis of an unaccelerated method \algname{D-EC-SGD} with Bias Correction and Double Contractive Compression~\citep{gao2024econtrol}. It is unknown if this is merely a defficiency in the analysis, or some inherent limitation of the method. The \algname{EControl} method resolved this issue in the unaccelerated case using a delicate error-control mechanism, but it seems to be difficult to adapt to cases where stepsizes are changing drastically which is precisely the case with acceleration. We leave it for future work to investigate if more refined analysis that improves the $\delta$ dependence in the deterministic term for \algname{ADEF} is possible or if the \algname{EControl} mechanism can be adapted to the accelerated setting.
\end{remark}

\begin{figure*}[t]
    \centering
    \includegraphics[width=0.34\textwidth]{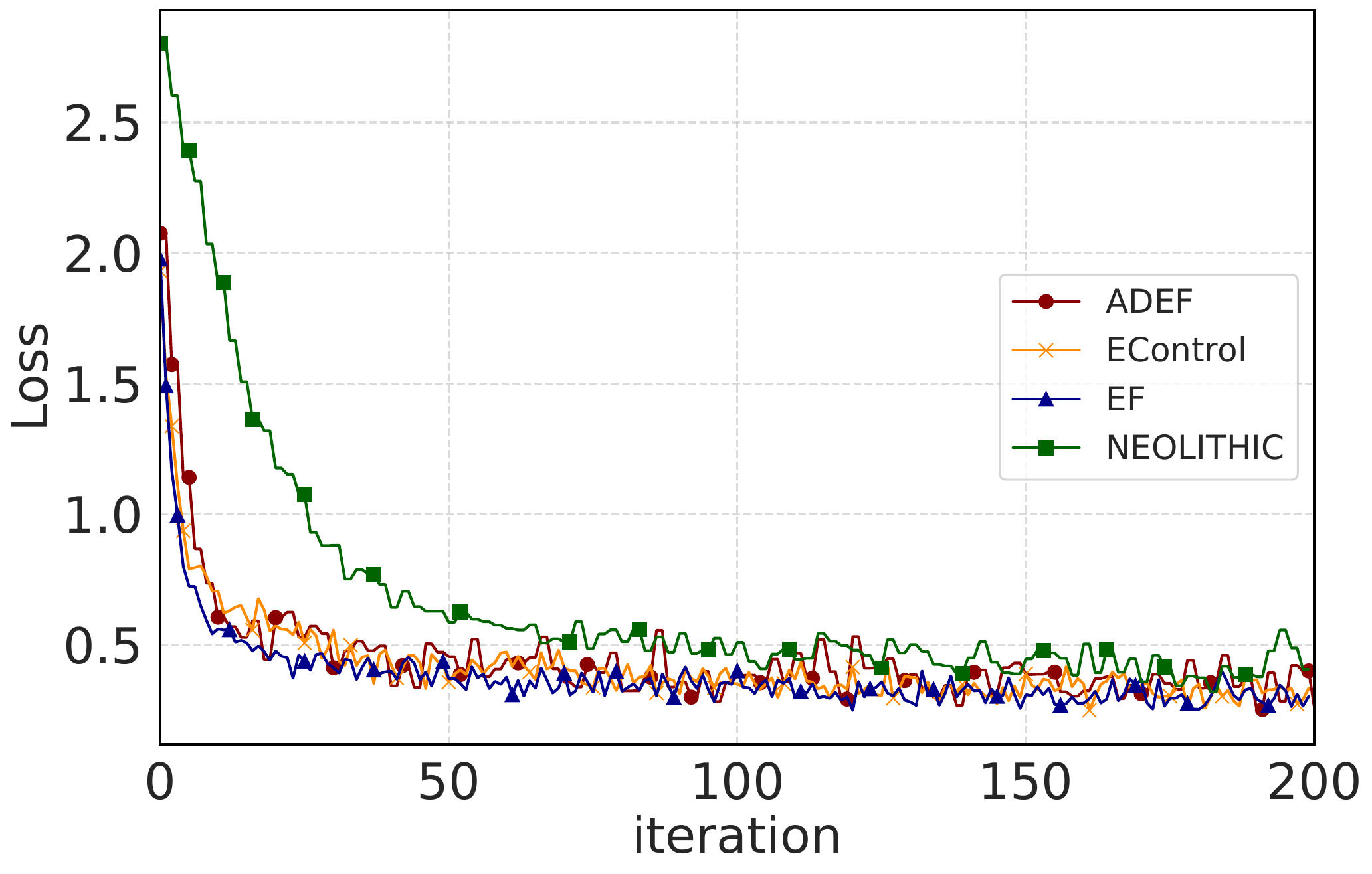}
    \hspace{4em}
    \includegraphics[width=0.34\textwidth]{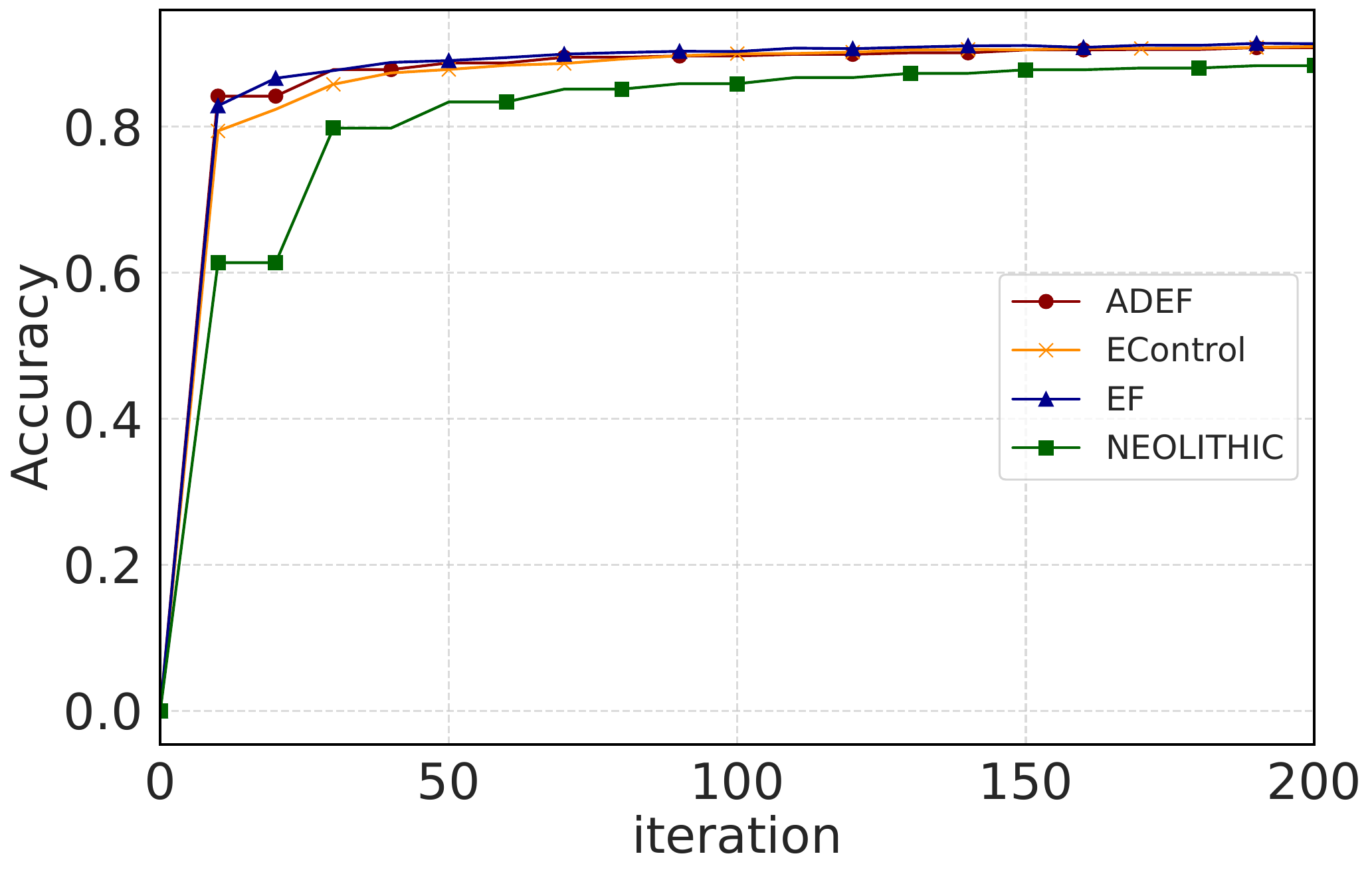}
    \vskip-5pt
    \caption{\textbf{Competitive performance.} %
         Comparison of the performance of \algname{ADEF}, \algname{EConrol},\algname{EF} and \algname{NEOLITHIC} on the MNIST classification problem. We use Top-$K$ compression with $\delta=0.1$. We see that \algname{ADEF} performs competetively in both the loss and the accuracy, while \algname{NEOLITHIC} performs the worst.}
    \label{fig:mnist-classification}
\end{figure*}

\section{Experiments}
\label{sec:experiments}

In this section we present numerical evaluations of our algorithm, and verify its theoretical properties discussed in the previous sections. We consider parameters of the following form: $a_t=\gamma(t+\frac{1}{\delta})$ and $A_0=\frac{1}{\delta^2}$, where $\gamma$ acts as the inverse of $M$ in \Cref{thm:convergence-rate}. For all experiments, we perform a grid search for $\gamma$, unless otherwise stated. We always report the performance in terms of the number of communications performed.

\subsection{Sythetic Logistic Regression}
\label{sec:synthetic-logistic-regression}

We first evaluate the performance of our algorithm on a synthetic logistic regression. We generate non-separable data with sklearn's built-in function~\citep{scikit-learn} and split the dataset into $n$ clients. We consider the following logistic loss for each data point $(\aa_i,b_i)$:
\begin{equation}
    \label{eq:logistic-loss}
    \ell((\aa_i,b_i),\xx) \eqdef -b_i\aa_i^\top\xx+\log(1+\exp(\aa_i^\top\xx)).
\end{equation}
For the stochastic oracle, we generate the noise from a Normal distribution with variance $\sigma^2$. We use the Top-$K$ compressor, with $\delta=0.1$. Now we verify some of the key properties of our algorithm.

\begin{figure}[tb]
    \centering
    \includegraphics[width=0.35\textwidth]{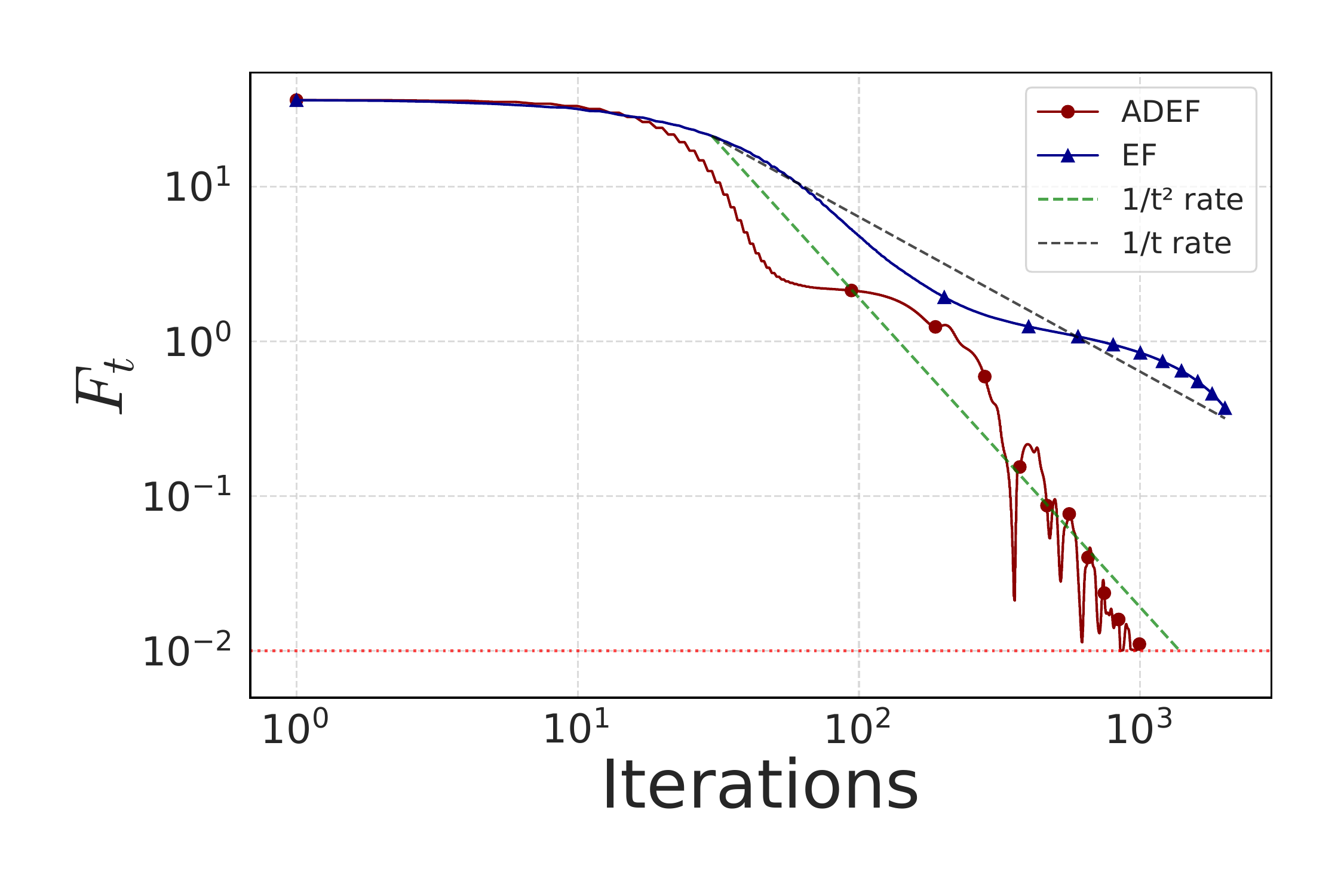}
    \caption{\textbf{Achieving acceleration.} Performance of \algname{ADEF} and \algname{EF} with $\sigma^2=0$. We set a target error of $0.01$. We see that \algname{ADEF} achieves the accelerated $\cO(\nicefrac{1}{T^2})$ rate.}
    \label{fig:acc}
\end{figure}

\textbf{Accelerated rate in the deterministic case} We first demonstrate that \Cref{alg:dist-acc-ef} indeed achieves the accelerated $\cO\left(\frac{1}{T^2}\right)$ rate when $\sigma^2=0$. We constrast this against the convergence rate of the unaccelerated \algname{EF}. The results are summarized in \Cref{fig:acc}. We see that \algname{ADEF} indeed achieves the accelerated rate $\cO\left(\frac{1}{T^2}\right)$, while \algname{EF} converges at the rate of $\cO\left(\frac{1}{T}\right)$.

\begin{figure}[tb]
    \centering
    \includegraphics[width=0.35\textwidth]{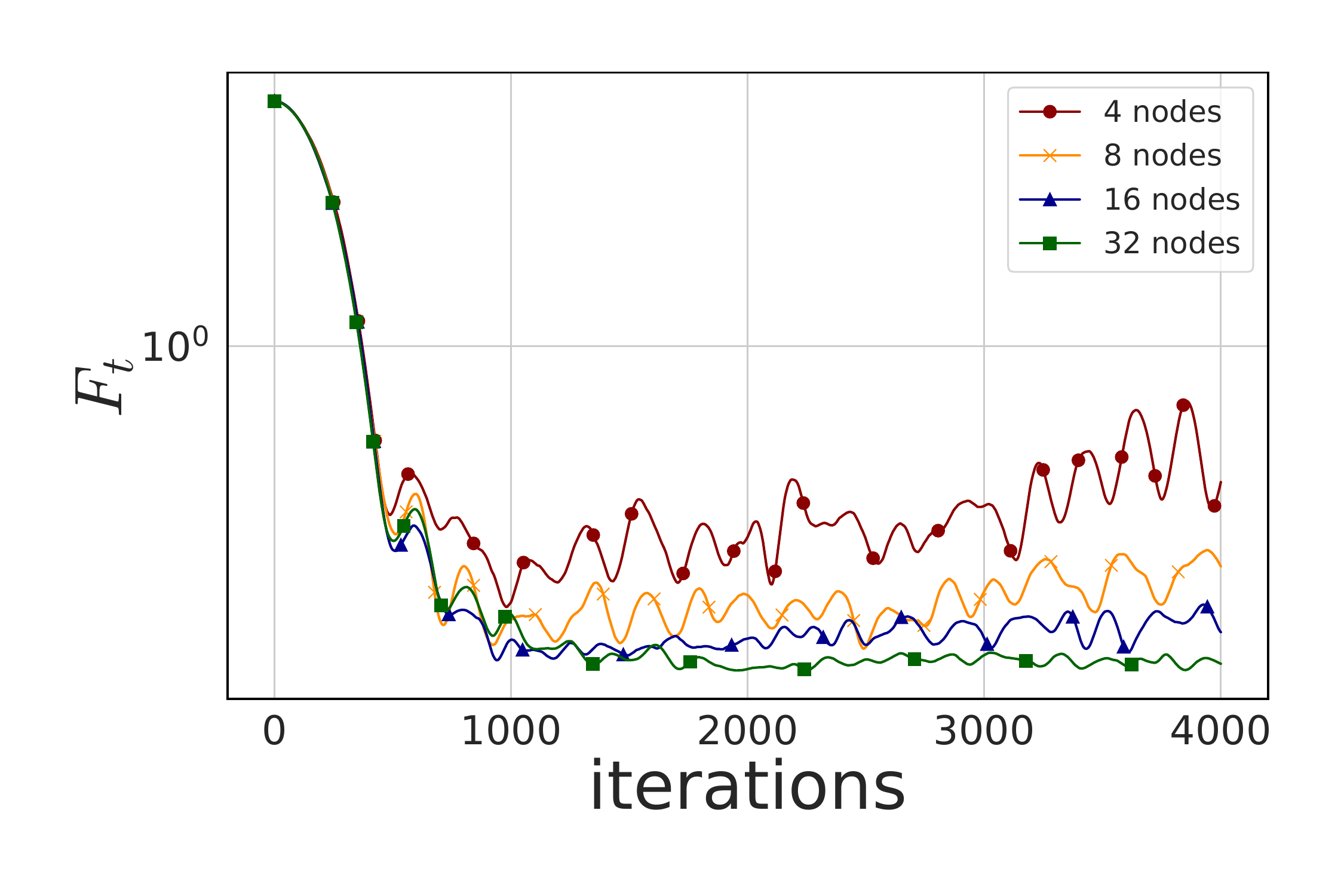}
    \caption{\textbf{Achieving linear speedup.} The performance of \algname{ADEF} with increasing number of clients $n$ for the synthetic logistic regressoin problem. We fix $\gamma$ to be $0.0001$. The error that the algorithm stabalizes around decreases as $n$ increases.}
    \label{fig:linear-speedup}
\end{figure}

\textbf{Linear speedup with the number of clients} Next we verify the linear speedup property of our algorithm with the number of clients $n$ with the leading terms in the convergence rate. We fix $\sigma^2=25$ and small enough $\gamma$ to be $0.0001$, and increase the number of clients $n$. The results are summarized in \Cref{fig:linear-speedup}. We see that the algorithm gets more stable as the number of clients increases, and the error that the algorithm stabalizes around decreases as $n$ increases, which verifies the linear speedup property of our algorithm.

\subsection{MNIST Classification}
\label{sec:mnist-classification}
Next we evaluate our algorithm on the MNIST classification problem~\citep{deng2012mnist} when using Top-$K$ compression with $\delta=0.1$, and compare it against the state of the art distributed optimization algorithms with communication compression, in particular, the classical \algname{EF} algorithm, and the recent \algname{EControl} algorithm. We also compare against the \algname{NEOLITHIC} method, where we set the number of repeated communication rounds to be $2$. The results are summarized in \Cref{fig:mnist-classification}. We see that \algname{ADEF} (with $\gamma=0.05$ as selected by grid search) performs competetively in both the loss and the accuracy, while \algname{NEOLITHIC} performs the worst.

\vspace{-0.5em}
\section{Conclusion and Outlook}
\label{sec:conclusion}

In this work, we propose a novel method \algname{ADEF} that provably achieves the accelerated convergence rate in the general convex regime for distributed optimization with contractive compression, making a first step towards the resolution of the open theoretical problem therein. We provide a general analysis framework for accelerated methods with inexact updates and apply it to analyze \algname{ADEF}, along with several others. Future works can look into the compression error dependence of algorithms employing double compressions such as ours; or study the possibility of extending our method, or indeed, (unaccelerated) error feedback methods, to constraint and composite optimization problems, which has been a long-standing theoretical issue in the field.

\bibliography{reference}

\begin{thebibliography}{60}
\providecommand{\natexlab}[1]{#1}
\providecommand{\url}[1]{\texttt{#1}}
\expandafter\ifx\csname urlstyle\endcsname\relax
  \providecommand{\doi}[1]{doi: #1}\else
  \providecommand{\doi}{doi: \begingroup \urlstyle{rm}\Url}\fi

\bibitem[Albasyoni et~al.(2020)Albasyoni, Safaryan, Condat, and Richt{\'a}rik]{albasyoni2020optimal}
Alyazeed Albasyoni, Mher Safaryan, Laurent Condat, and Peter Richt{\'a}rik.
\newblock Optimal gradient compression for distributed and federated learning.
\newblock \emph{arXiv preprint arXiv:2010.03246}, 2020.

\bibitem[Alistarh et~al.(2017)Alistarh, Grubic, Li, Tomioka, and Vojnovic]{QSGD}
Dan Alistarh, Demjan Grubic, Jerry Li, Ryota Tomioka, and Milan Vojnovic.
\newblock {QSGD}: Communication-efficient {SGD} via gradient quantization and encoding.
\newblock In \emph{Proceedings of Advances in Neural Information Processing Systems}, pages 1709--1720, 2017.

\bibitem[Alistarh et~al.(2018)Alistarh, Hoefler, Johansson, Konstantinov, Khirirat, and Renggli]{alistarh2018convergence}
Dan Alistarh, Torsten Hoefler, Mikael Johansson, Nikola Konstantinov, Sarit Khirirat, and C{\'e}dric Renggli.
\newblock The convergence of sparsified gradient methods.
\newblock In \emph{Proceedings of Advances in Neural Information Processing Systems}, 2018.

\bibitem[Allen-Zhu(2018)]{allen2018katyusha}
Zeyuan Allen-Zhu.
\newblock Katyusha: The first direct acceleration of stochastic gradient methods.
\newblock \emph{Journal of Machine Learning Research}, 18\penalty0 (221):\penalty0 1--51, 2018.

\bibitem[Bernstein et~al.(2018)Bernstein, Wang, Azizzadenesheli, and Anandkumar]{bernstein2018signsgd}
Jeremy Bernstein, Yu-Xiang Wang, Kamyar Azizzadenesheli, and Animashree Anandkumar.
\newblock signsgd: Compressed optimisation for non-convex problems.
\newblock In \emph{International Conference on Machine Learning}, pages 560--569. PMLR, 2018.

\bibitem[Beznosikov et~al.(2020)Beznosikov, Horv{\'a}th, Richt{\'a}rik, and Safaryan]{beznosikov2020biased}
Aleksandr Beznosikov, Samuel Horv{\'a}th, Peter Richt{\'a}rik, and Mher Safaryan.
\newblock On biased compression for distributed learning.
\newblock \emph{Journal on Machine Learning Research}, 2020.

\bibitem[Beznosikov et~al.(2023)Beznosikov, Horv{\'a}th, Richt{\'a}rik, and Safaryan]{beznosikov2023biased}
Aleksandr Beznosikov, Samuel Horv{\'a}th, Peter Richt{\'a}rik, and Mher Safaryan.
\newblock On biased compression for distributed learning.
\newblock \emph{Journal of Machine Learning Research}, 24\penalty0 (276):\penalty0 1--50, 2023.

\bibitem[Bylinkin and Beznosikov(2024)]{bylinkin2024accelerated}
Dmitry Bylinkin and Aleksandr Beznosikov.
\newblock Accelerated methods with compressed communications for distributed optimization problems under data similarity.
\newblock \emph{arXiv preprint arXiv:2412.16414}, 2024.

\bibitem[Cordonnier(2018)]{cordonnier2018convex}
Jean-Baptiste Cordonnier.
\newblock Convex optimization using sparsified stochastic gradient descent with memory.
\newblock Technical report, 2018.

\bibitem[Deng(2012)]{deng2012mnist}
Li~Deng.
\newblock The mnist database of handwritten digit images for machine learning research.
\newblock \emph{IEEE Signal Processing Magazine}, 2012.

\bibitem[Devolder et~al.(2014)Devolder, Glineur, and Nesterov]{devolder2014first}
Olivier Devolder, Fran{\c{c}}ois Glineur, and Yurii Nesterov.
\newblock First-order methods of smooth convex optimization with inexact oracle.
\newblock \emph{Mathematical Programming}, 146:\penalty0 37--75, 2014.

\bibitem[Fatkhullin et~al.(2021)Fatkhullin, Sokolov, Gorbunov, Li, and Richtárik]{fatkhullin2021ef21}
Ilyas Fatkhullin, Igor Sokolov, Eduard Gorbunov, Zhize Li, and Peter Richtárik.
\newblock Ef21 with bells \& whistles: Practical algorithmic extensions of modern error feedback.
\newblock \emph{arXiv preprint arXiv: 2110.03294}, 2021.

\bibitem[Fatkhullin et~al.(2023)Fatkhullin, Tyurin, and Richt{\'a}rik]{fatkhullin2023momentum}
Ilyas Fatkhullin, Alexander Tyurin, and Peter Richt{\'a}rik.
\newblock Momentum provably improves error feedback!
\newblock \emph{arXiv preprint arXiv: 2305.15155}, 2023.

\bibitem[Gao et~al.(2024)Gao, Islamov, and Stich]{gao2024econtrol}
Yuan Gao, Rustem Islamov, and Sebastian~U Stich.
\newblock {EC}ontrol: Fast distributed optimization with compression and error control.
\newblock In \emph{The Twelfth International Conference on Learning Representations}, 2024.

\bibitem[Gorbunov et~al.(2020)Gorbunov, Kovalev, Makarenko, and Richt{\'a}rik]{gorbunov2020linearly}
Eduard Gorbunov, Dmitry Kovalev, Dmitry Makarenko, and Peter Richt{\'a}rik.
\newblock Linearly converging error compensated sgd.
\newblock In \emph{Proceedings of Advances in Neural Information Processing Systems}, 2020.

\bibitem[Gupta et~al.(2024)Gupta, Siegel, and Wojtowytsch]{gupta2024nesterov}
Kanan Gupta, Jonathan~W. Siegel, and Stephan Wojtowytsch.
\newblock Nesterov acceleration despite very noisy gradients.
\newblock In \emph{The Thirty-eighth Annual Conference on Neural Information Processing Systems}, 2024.

\bibitem[Gupta et~al.(2015)Gupta, Agrawal, Gopalakrishnan, and Narayanan]{gupta2015deep}
Suyog Gupta, Ankur Agrawal, Kailash Gopalakrishnan, and Pritish Narayanan.
\newblock Deep learning with limited numerical precision.
\newblock In \emph{International conference on machine learning}, pages 1737--1746. PMLR, 2015.

\bibitem[He et~al.(2023)He, Huang, Chen, Yin, and Yuan]{he2023lower}
Yutong He, Xinmeng Huang, Yiming Chen, Wotao Yin, and Kun Yuan.
\newblock Lower bounds and accelerated algorithms in distributed stochastic optimization with communication compression.
\newblock \emph{arXiv preprint arXiv: 2305.07612}, 2023.

\bibitem[Horváth et~al.(2019)Horváth, Ho, Horváth, Sahu, Canini, and Richtárik]{horvath2019natural}
Samuel Horváth, Chen-Yu Ho, \v{L}udovít Horváth, Atal~Narayan Sahu, Marco Canini, and Peter Richtárik.
\newblock Natural compression for distributed deep learning.
\newblock \emph{arXiv preprint arXiv: 1905.10988}, 2019.

\bibitem[Huang et~al.(2022)Huang, Chen, Yin, and Yuan]{huang2022lowerbounds}
Xinmeng Huang, Yiming Chen, Wotao Yin, and Kun Yuan.
\newblock Lower bounds and nearly optimal algorithms in distributed learning with communication compression.
\newblock In \emph{Proceedings of Advances in Neural Information Processing Systems}, 2022.

\bibitem[Islamov et~al.(2023)Islamov, Qian, Hanzely, Safaryan, and Richt{\'a}rik]{islamov2023distributed}
Rustem Islamov, Xun Qian, Slavom{\'\i}r Hanzely, Mher Safaryan, and Peter Richt{\'a}rik.
\newblock Distributed newton-type methods with communication compression and bernoulli aggregation.
\newblock \emph{Transactions on Machine Learning Research}, 2023.

\bibitem[Jiang et~al.(2024)Jiang, Rodomanov, and Stich]{jiang2024federated}
Xiaowen Jiang, Anton Rodomanov, and Sebastian~U Stich.
\newblock Federated optimization with doubly regularized drift correction.
\newblock \emph{arXiv preprint arXiv:2404.08447}, 2024.

\bibitem[Kairouz et~al.(2019)Kairouz, McMahan, Avent, Bellet, Bennis, Bhagoji, Bonawitz, Charles, Cormode, Cummings, D'Oliveira, Rouayheb, Evans, Gardner, Garrett, Gasc{\'{o}}n, Ghazi, Gibbons, Gruteser, Harchaoui, He, He, Huo, Hutchinson, Hsu, Jaggi, Javidi, Joshi, Khodak, Kone{\v{c}}n{\'y}, Korolova, Koushanfar, Koyejo, Lepoint, Liu, Mittal, Mohri, Nock, {\"{O}}zg{\"{u}}r, Pagh, Raykova, Qi, Ramage, Raskar, Song, Song, Stich, Sun, Suresh, Tram{\`{e}}r, Vepakomma, Wang, Xiong, Xu, Yang, Yu, Yu, and Zhao]{kairouz2019advances}
Peter Kairouz, H.~Brendan McMahan, Brendan Avent, Aur{\'{e}}lien Bellet, Mehdi Bennis, Arjun~Nitin Bhagoji, Kallista~A. Bonawitz, Zachary Charles, Graham Cormode, Rachel Cummings, Rafael G.~L. D'Oliveira, Salim~El Rouayheb, David Evans, Josh Gardner, Zachary Garrett, Adri{\`{a}} Gasc{\'{o}}n, Badih Ghazi, Phillip~B. Gibbons, Marco Gruteser, Za{\"{\i}}d Harchaoui, Chaoyang He, Lie He, Zhouyuan Huo, Ben Hutchinson, Justin Hsu, Martin Jaggi, Tara Javidi, Gauri Joshi, Mikhail Khodak, Jakub Kone{\v{c}}n{\'y}, Aleksandra Korolova, Farinaz Koushanfar, Sanmi Koyejo, Tancr{\`{e}}de Lepoint, Yang Liu, Prateek Mittal, Mehryar Mohri, Richard Nock, Ayfer {\"{O}}zg{\"{u}}r, Rasmus Pagh, Mariana Raykova, Hang Qi, Daniel Ramage, Ramesh Raskar, Dawn Song, Weikang Song, Sebastian~U. Stich, Ziteng Sun, Ananda~Theertha Suresh, Florian Tram{\`{e}}r, Praneeth Vepakomma, Jianyu Wang, Li~Xiong, Zheng Xu, Qiang Yang, Felix~X. Yu, Han Yu, and Sen Zhao.
\newblock Advances and open problems in federated learning.
\newblock \emph{arXiv preprint arXiv: 1912.04977}, 2019.

\bibitem[Karimireddy et~al.(2019)Karimireddy, Rebjock, Stich, and Jaggi]{karimireddy2019error}
Sai~Praneeth Karimireddy, Quentin Rebjock, Sebastian Stich, and Martin Jaggi.
\newblock Error feedback fixes signsgd and other gradient compression schemes.
\newblock In \emph{Proceedings of the 36th International Conference on Machine Learning (ICML 2019)}, 2019.

\bibitem[Khaled and Jin(2022)]{khaled2022faster}
Ahmed Khaled and Chi Jin.
\newblock Faster federated optimization under second-order similarity.
\newblock \emph{arXiv preprint arXiv:2209.02257}, 2022.

\bibitem[Koloskova et~al.(2019)Koloskova, Stich, and Jaggi]{koloskova2019Decentralized}
Anastasia Koloskova, Sebastian Stich, and Martin Jaggi.
\newblock Decentralized stochastic optimization and gossip algorithms with compressed communication.
\newblock In \emph{Proceedings of the 36th International Conference on Machine Learning}, 2019.

\bibitem[Koloskova et~al.(2020)Koloskova, Lin, Stich, and Jaggi]{koloskova2020Decentralized}
Anastasia Koloskova, Tao Lin, Sebastian Stich, and Martin Jaggi.
\newblock Decentralized deep learning with arbitrary communication compression.
\newblock In \emph{Proceedings of International Conference on Learning Representations}, 2020.

\bibitem[Kone\v{c}n\'{y} et~al.(2016)Kone\v{c}n\'{y}, McMahan, Yu, Richt\'{a}rik, Suresh, and Bacon]{FEDLEARN}
Jakub Kone\v{c}n\'{y}, H.~Brendan McMahan, Felix Yu, Peter Richt\'{a}rik, Ananda~Theertha Suresh, and Dave Bacon.
\newblock Federated learning: strategies for improving communication efficiency.
\newblock In \emph{Proceedings of NIPS Private Multi-Party Machine Learning Workshop}, 2016.

\bibitem[Li and Li(2023)]{li2023analysis}
Xiaoyun Li and Ping Li.
\newblock Analysis of error feedback in federated non-convex optimization with biased compression: Fast convergence and partial participation.
\newblock In \emph{Proceedings of the 40th International Conference on Machine Learning}, 2023.

\bibitem[Li and Richtárik(2021)]{li2021canita}
Zhize Li and Peter Richtárik.
\newblock Canita: Faster rates for distributed convex optimization with communication compression, 2021.

\bibitem[Li et~al.(2020)Li, Kovalev, Qian, and Richtárik]{li2020acceleration}
Zhize Li, Dmitry Kovalev, Xun Qian, and Peter Richtárik.
\newblock Acceleration for compressed gradient descent in distributed and federated optimization, 2020.

\bibitem[Lian et~al.(2017)Lian, Zhang, Zhang, Hsieh, Zhang, and Liu]{lian2017can}
Xiangru Lian, Ce~Zhang, Huan Zhang, Cho-Jui Hsieh, Wei Zhang, and Ji~Liu.
\newblock Can decentralized algorithms outperform centralized algorithms? a case study for decentralized parallel stochastic gradient descent.
\newblock In \emph{Proceedings of Advances in Neural Information Processing Systems (NIPS)}, 2017.

\bibitem[Lin et~al.(2018)Lin, Han, Mao, Wang, and Dally]{lin2018deep}
Yujun Lin, Song Han, Huizi Mao, Yu~Wang, and Bill Dally.
\newblock Deep gradient compression: Reducing the communication bandwidth for distributed training.
\newblock In \emph{International Conference on Learning Representations}, 2018.

\bibitem[Mishchenko et~al.(2019)Mishchenko, Gorbunov, Tak{\'a}{\v{c}}, and Richt{\'a}rik]{DIANA}
Konstantin Mishchenko, Eduard Gorbunov, Martin Tak{\'a}{\v{c}}, and Peter Richt{\'a}rik.
\newblock Distributed learning with compressed gradient differences.
\newblock \emph{arXiv preprint arXiv:1901.09269}, 2019.

\bibitem[Murata and Suzuki(2019)]{murata2019accelerated}
Tomoya Murata and Taiji Suzuki.
\newblock Accelerated sparsified sgd with error feedback.
\newblock \emph{arXiv preprint arXiv:1905.12224}, 2019.

\bibitem[Nesterov(1983)]{nesterov1983method}
Yurii Nesterov.
\newblock A method for solving the convex programming problem with convergence rate o (1/k2).
\newblock In \emph{Dokl akad nauk Sssr}, volume 269, page 543, 1983.

\bibitem[Nesterov et~al.(2018)]{nesterov2018lectures}
Yurii Nesterov et~al.
\newblock \emph{Lectures on convex optimization}, volume 137.
\newblock Springer, 2018.

\bibitem[Paszke et~al.(2019)Paszke, Gross, Massa, Lerer, Bradbury, Chanan, Killeen, Lin, Gimelshein, Antiga, Desmaison, Kopf, Yang, DeVito, Raison, Tejani, Chilamkurthy, Steiner, Fang, Bai, and Chintala]{pazske2019pytorch}
Adam Paszke, Sam Gross, Francisco Massa, Adam Lerer, James Bradbury, Gregory Chanan, Trevor Killeen, Zeming Lin, Natalia Gimelshein, Luca Antiga, Alban Desmaison, Andreas Kopf, Edward Yang, Zachary DeVito, Martin Raison, Alykhan Tejani, Sasank Chilamkurthy, Benoit Steiner, Lu~Fang, Junjie Bai, and Soumith Chintala.
\newblock Pytorch: An imperative style, high-performance deep learning library.
\newblock In \emph{Proceedings of Advances in Neural Information Processing Systems 32}, 2019.

\bibitem[Pedregosa et~al.(2011)Pedregosa, Varoquaux, Gramfort, Michel, Thirion, Grisel, Blondel, Prettenhofer, Weiss, Dubourg, Vanderplas, Passos, Cournapeau, Brucher, Perrot, and Duchesnay]{scikit-learn}
F.~Pedregosa, G.~Varoquaux, A.~Gramfort, V.~Michel, B.~Thirion, O.~Grisel, M.~Blondel, P.~Prettenhofer, R.~Weiss, V.~Dubourg, J.~Vanderplas, A.~Passos, D.~Cournapeau, M.~Brucher, M.~Perrot, and E.~Duchesnay.
\newblock Scikit-learn: Machine learning in {P}ython.
\newblock \emph{Journal of Machine Learning Research}, 12:\penalty0 2825--2830, 2011.

\bibitem[Qian et~al.(2021{\natexlab{a}})Qian, Richt{\'a}rik, and Zhang]{qian2021error}
Xun Qian, Peter Richt{\'a}rik, and Tong Zhang.
\newblock Error compensated distributed sgd can be accelerated.
\newblock \emph{Advances in Neural Information Processing Systems}, 34:\penalty0 30401--30413, 2021{\natexlab{a}}.

\bibitem[Qian et~al.(2021{\natexlab{b}})Qian, Richt{\'a}rik, and Zhang]{qian2021errorcompensatedSGD}
Xun Qian, Peter Richt{\'a}rik, and Tong Zhang.
\newblock Error compensated distributed sgd can be accelerated.
\newblock In \emph{Proceedings of Advances in Neural Information Processing Systems}, 2021{\natexlab{b}}.

\bibitem[Ramesh et~al.(2021)Ramesh, Pavlov, Goh, Gray, Voss, Radford, Chen, and Sutskever]{ramesh2021zero}
Aditya Ramesh, Mikhail Pavlov, Gabriel Goh, Scott Gray, Chelsea Voss, Alec Radford, Mark Chen, and Ilya Sutskever.
\newblock Zero-shot text-to-image generation.
\newblock In \emph{International Conference on Machine Learning}, pages 8821--8831. PMLR, 2021.

\bibitem[Ramesh et~al.(2022)Ramesh, Dhariwal, Nichol, Chu, and Chen]{ramesh2022hierarchical}
Aditya Ramesh, Prafulla Dhariwal, Alex Nichol, Casey Chu, and Mark Chen.
\newblock Hierarchical text-conditional image generation with clip latents.
\newblock 2022.

\bibitem[Richt{\'a}rik et~al.(2021)Richt{\'a}rik, Sokolov, and Fatkhullin]{richtarik2021ef21}
Peter Richt{\'a}rik, Igor Sokolov, and Ilyas Fatkhullin.
\newblock Ef21: A new, simpler, theoretically better, and practically faster error feedback.
\newblock In \emph{Proceedings of Advances in Neural Information Processing Systems}, 2021.

\bibitem[Rodomanov et~al.(2024)Rodomanov, Jiang, and Stich]{rodomanov2024universality}
Anton Rodomanov, Xiaowen Jiang, and Sebastian~U Stich.
\newblock Universality of adagrad stepsizes for stochastic optimization: Inexact oracle, acceleration and variance reduction.
\newblock In \emph{The Thirty-eighth Annual Conference on Neural Information Processing Systems}, 2024.

\bibitem[Safaryan et~al.(2022)Safaryan, Islamov, Qian, and Richtarik]{safaryanFedNL2022}
Mher Safaryan, Rustem Islamov, Xun Qian, and Peter Richtarik.
\newblock {F}ed{NL}: Making {N}ewton-type methods applicable to federated learning.
\newblock In \emph{Proceedings of the 39th International Conference on Machine Learning}, 2022.

\bibitem[Sahu et~al.(2021)Sahu, Dutta, M~Abdelmoniem, Banerjee, Canini, and Kalnis]{sahu2021rethinking}
Atal Sahu, Aritra Dutta, Ahmed M~Abdelmoniem, Trambak Banerjee, Marco Canini, and Panos Kalnis.
\newblock Rethinking gradient sparsification as total error minimization.
\newblock \emph{Advances in Neural Information Processing Systems}, 34:\penalty0 8133--8146, 2021.

\bibitem[Sapio et~al.(2021)Sapio, Canini, Ho, Nelson, Kalnis, Kim, Krishnamurthy, Moshref, Ports, and Richt{\'a}rik]{sapio2021scaling}
Amedeo Sapio, Marco Canini, Chen-Yu Ho, Jacob Nelson, Panos Kalnis, Changhoon Kim, Arvind Krishnamurthy, Masoud Moshref, Dan Ports, and Peter Richt{\'a}rik.
\newblock Scaling distributed machine learning with $\{$In-Network$\}$ aggregation.
\newblock In \emph{18th USENIX Symposium on Networked Systems Design and Implementation (NSDI 21)}, pages 785--808, 2021.

\bibitem[Seide et~al.(2014)Seide, Fu, Droppo, Li, and Yu]{seide20141}
Frank Seide, Hao Fu, Jasha Droppo, Gang Li, and Dong Yu.
\newblock 1-bit stochastic gradient descent and its application to data-parallel distributed training of speech dnns.
\newblock In \emph{Proceedings of 15th annual conference of the international speech communication association}, 2014.

\bibitem[Shoeybi et~al.(2019)Shoeybi, Patwary, Puri, Gresley, Casper, and Catanzaro]{shoeybi2019megatron}
Mohammad Shoeybi, Mostofa Patwary, Raul Puri, Patrick~Le Gresley, Jared Casper, and Bryan Catanzaro.
\newblock Megatron-lm: Training multi-billion parameter language models using model parallelism.
\newblock \emph{arXiv preprint arXiv:1909.08053}, 2019.

\bibitem[Stich(2020)]{stich2020communication}
Sebastian~U. Stich.
\newblock On communication compression for distributed optimization on heterogeneous data.
\newblock \emph{arXiv preprint arXiv: 2009.02388}, 2020.

\bibitem[Stich and Karimireddy(2020)]{stich2020error}
Sebastian~U Stich and Sai~Praneeth Karimireddy.
\newblock The error-feedback framework: Better rates for sgd with delayed gradients and compressed updates.
\newblock \emph{Journal of Machine Learning Research}, 2020.

\bibitem[Stich et~al.(2018)Stich, Cordonnier, and Jaggi]{stich2018sparsified}
Sebastian~U Stich, Jean-Baptiste Cordonnier, and Martin Jaggi.
\newblock Sparsified sgd with memory.
\newblock In \emph{Proceedings of Advances in Neural Information Processing Systems}, 2018.

\bibitem[Strom(2015)]{strom15_interspeech}
Nikko Strom.
\newblock {Scalable distributed DNN training using commodity GPU cloud computing}.
\newblock In \emph{Proceedings of Interspeech 2015}, 2015.

\bibitem[Sun et~al.(2019)Sun, Shao, Jiang, Cui, Lei, Xu, and Wang]{sun2019sparse}
Haobo Sun, Yingxia Shao, Jiawei Jiang, Bin Cui, Kai Lei, Yu~Xu, and Jiang Wang.
\newblock Sparse gradient compression for distributed sgd.
\newblock In \emph{International Conference on Database Systems for Advanced Applications}, pages 139--155. Springer, 2019.

\bibitem[Vogels et~al.(2019)Vogels, Karimireddy, and Jaggi]{vogels2019PowerSGD}
Thijs Vogels, Sai~Praneeth Karimireddy, and Martin Jaggi.
\newblock Power{SGD}: Practical low-rank gradient compression for distributed optimization.
\newblock In \emph{Advances in Neural Information Processing Systems 32}, 2019.

\bibitem[Wang et~al.(2020)Wang, Fu, He, Hao, and Wu]{wang2020survey}
Meng Wang, Weijie Fu, Xiangnan He, Shijie Hao, and Xindong Wu.
\newblock A survey on large-scale machine learning.
\newblock \emph{IEEE Transactions on Knowledge and Data Engineering}, 2020.

\bibitem[Wen et~al.(2017)Wen, Xu, Yan, Wu, Wang, Chen, and Li]{wen2017terngrad}
Wei Wen, Cong Xu, Feng Yan, Chunpeng Wu, Yandan Wang, Yiran Chen, and Hai Li.
\newblock Terngrad: Ternary gradients to reduce communication in distributed deep learning.
\newblock \emph{Advances in neural information processing systems}, 30, 2017.

\bibitem[Zheng et~al.(2019)Zheng, Huang, and Kwok]{zheng2019communication}
Shuai Zheng, Ziyue Huang, and James Kwok.
\newblock Communication-efficient distributed blockwise momentum sgd with error-feedback.
\newblock \emph{Advances in Neural Information Processing Systems}, 32, 2019.

\bibitem[Zhu et~al.(2023)Zhu, Ghosh, and Mazumdar]{zhu2023optimal}
Heng Zhu, Avishek Ghosh, and Arya Mazumdar.
\newblock Optimal compression of unit norm vectors in the high distortion regime.
\newblock In \emph{2023 IEEE International Symposium on Information Theory (ISIT)}, pages 719--724. IEEE, 2023.

\end{thebibliography}
\bibliographystyle{plainnat}

\appendix
\onecolumn

\section{Auxilary Facts and Results}
\label{sec:auxilary}
In this section we collect some auxilary facts and results that are useful for the analysis of our algorithms. The first one is a simple fact regarding the square of the norm of a sum of vectors.

\begin{fact}
    \label{fact:square-sum}
    For any $\gamma_1,\dots,\gamma_T$, we have:
    \begin{equation}
        \label{eq:square-sum}
        \norm{\sum_{t=1}^T\gamma_t}^2\leq T\sum_{t=1}^T\norm{\gamma_t}^2.
    \end{equation}
\end{fact}

The next one is on upper bounding the sum of a sequence under a descent condition on the sequence:

\begin{restatable}{lemma}{sumLemma}
    \label{lem:sum-lemma}
    Given a sequence of non-negative values $\{\alpha_t\}_{t\in[T-1]}$, and some other sequences $\{\lambda_t\}_{t\in[T-1]}$. If there exists $\gamma\in(0,1)$ such that the following holds:
    \begin{equation}
        \label{eq:sum-lemma-assumption}
        \alpha_{t+1}\leq (1-\gamma)\alpha_t + \lambda_t, \quad \alpha_0=0,\,
    \end{equation}
    then we have:
    \begin{equation}
        \label{eq:sum-lemma}
        \sum_{t=0}^{T-1}\alpha_{t+1}\leq \frac{1}{\gamma}\sum_{t=0}^{T-1}\lambda_t
    \end{equation}
\end{restatable}

\begin{proof}
    We sum \Cref{eq:sum-lemma-assumption} on both sides from $t=0$ to $T-1$:
    \begin{align*}
        \sum_{t=0}^{T-1}\alpha_{t+1} &\leq (1-\gamma)\sum_{t=0}^{T-1}\alpha_t + \sum_{t=0}^{T-1}\lambda_t\\
        &\overset{(i)}{=} (1-\gamma)\sum_{t=1}^{T-1}\alpha_t + \sum_{t=0}^{T-1}\lambda_t,\,
    \end{align*}
    where in $(i)$ we used the fact that $\alpha_0=0$. Now we can subtract $(1-\gamma)\sum_{t=1}^{T-1}\alpha_t$ from both sides and get:
    \[
        \gamma \sum_{t=0}^{T-1}\alpha_{t+1} \leq \sum_{t=0}^{T-1}\lambda_t,\,
    \]
    divide both sides by $\gamma$ and we get the desired result.
\end{proof}
The next fact is on the convexity and smoothness of a function:
\begin{fact}
    \label{fact:convexity-smoothness}
    Given \Cref{assumption:convexity,assumption:smoothness}, for any $\xx,\yy\in\R^d$, we have:
    \begin{equation}
        \label{eq:convexity-smoothness}
        \norm{\nabla f(\xx)-\nabla f(\yy)}^2\leq 2L\beta_f(\xx,\yy).
    \end{equation}
\end{fact}
We also have the following fact.
\begin{fact}
    \label{fact:ell-L}
    If \Cref{assumption:f-i-smoothness} holds with some $\ell\geq 0$, then \Cref{assumption:smoothness} holds with $L\leq \ell$.
\end{fact}
See, e.g. Theorem 2.1.5 in \citet{nesterov2018lectures} for a proof for both \Cref{fact:convexity-smoothness} and \Cref{fact:ell-L}. 
\section{Analysis of Accelerated Method with Inexact Updates}
\label{sec:acc-inexact-proof}

In this section we present the proof of \Cref{thm:main-descent}.

\mainDescent*

\begin{proof}
  Recall that
  \[
    \vv_{t + 1} = \vv_t - a_{t + 1} \hat{\gg}_t.
  \]

  We consider $\wtilde\vv_t,\wtilde\vv_0=\vv_0$ defined the in the following way:
  \[
      \wtilde\vv_{t+1} \eqdef \wtilde\vv_t - a_{t+1}\gg_t
  \]

  We also consider $\wtilde R_t\eqdef \EbNSq{\wtilde\vv_t -\xx}$.

  We start by giving a one step descent. Given any $\xx\in \R^d$ we have the following:
  \begin{align*}
      & \Eb{A_tf(\xx_t)+a_{t+1}f(\xx)+\frac{1}{2}\norm{\wtilde\vv_t-\xx}^2}\\
  \overset{(i)}{=} & \Eb{A_t(f(\yy_t)+\inp{\nabla f(\yy_t)}{\xx_t-\yy_t}+\beta_f(\yy_t,\xx_t))}\\
      &  + \Eb{ a_{t+1}(f(\yy_t)+\inp{\nabla f(\yy_t)}{\xx-\yy_t}+\beta_f(\yy_t,\xx))+\frac{1}{2}\norm{\wtilde\vv_t-\xx}^2}\\
  \overset{(ii)}{=} & \Eb{A_t(f(\yy_t)+\inp{\nabla f(\yy_t)}{\xx_t-\yy_t}+\beta_f(\yy_t,\xx_t))}\\
      & + \Eb{ a_{t+1}(f(\yy_t)+\inp{\gg_t }{\xx-\yy_t}+\beta_f(\yy_t,\xx))+\frac{1}{2}\norm{\wtilde\vv_t-\xx}^2}\\
  \overset{(iii)}{=} & \Eb{A_t(f(\yy_t)+\inp{\nabla f(\yy_t)}{\xx_t-\yy_t}+\beta_f(\yy_t,\xx_t)) }\\
      & + \Eb{a_{t+1}(f(\yy_t)+\inp{\gg_t }{\wtilde\vv_{t+1}-\yy_t}+\beta_f(\yy_t,\xx))+\frac{1}{2}\norm{\wtilde\vv_t-\wtilde\vv_{t+1}}^2 + \frac{1}{2}\norm{\wtilde\vv_{t+1}-\xx}^2}\\
  = & \Eb{A_t(f(\yy_t)+\inp{\nabla f(\yy_t)}{\xx_t-\yy_t}+\beta_f(\yy_t,\xx_t)) }\\
      &  +\Eb{ a_{t+1}(f(\yy_t)+\inp{\nabla f(\yy_t)}{\wtilde\vv_{t+1}-\yy_t}+\beta_f(\yy_t,\xx))+\frac{1}{2}\norm{\wtilde\vv_t-\wtilde\vv_{t+1}}^2 }\\
  & + \Eb{\frac{1}{2}\norm{\wtilde\vv_{t+1}-\xx}^2 + a_{t+1}\inp{\gg_t -\nabla f(\yy_t)}{\wtilde\vv_{t+1}-\wtilde\vv_t}} \\
  \overset{(iv)}{=} & \Eb{A_{t+1}\left(f(\yy_t) + \inp{\nabla f(\yy_t)}{\xx_{t+1}-\yy_t}\right) + A_t\beta_f(\yy_t,\xx_t) + a_{t+1}\beta_f(\yy_t,\xx) +\frac{1}{2}\norm{\wtilde\vv_t-\wtilde\vv_{t+1}}^2 }\\
  & + \Eb{\frac{1}{2}\norm{\wtilde\vv_{t+1}-\xx}^2 + a_{t+1}\inp{\gg_t -\nabla f(\yy_t)}{\wtilde\vv_{t+1}-\wtilde\vv_t} - a_{t+1}\inp{\nabla f(\yy_t)}{\ee_{t+1}}} \,
  \end{align*}
  where for $(i)$ we applied the convexity of $f$. For $(ii)$ we applied the independence and unbiasedness of $\gg_t $. For $(iii)$ we used the fact that
  \[
      \wtilde\vv_{t+1} = \argmin_{\vv\in\R^d} \Psi_t(\vv_t)\eqdef a_{t+1} \inp{\gg_t }{\vv} + \frac{1}{2}\norm{\vv-\wtilde\vv_t}^2
  \]
  where $\Psi_t$ is $1$-strongly convex. For $(iv)$ we used the definition of $\xx_{t+1}$. Next, we apply the smoothness and convexity of $f$:
  \begin{align*}
      f(\yy_t) + \inp{\nabla f(\yy_t)}{\xx_{t+1}-\yy_t} & = f(\xx_{t+1}) -\inp{\nabla f(\xx_{t+1})-\nabla f(\yy_t)}{\xx_{t+1}-\yy_t} + \beta_f(\xx_{t+1},\yy_t)\\
      &\geq f(\xx_{t+1}) + \beta_f(\xx_{t+1},\yy_t) - L\norm{\xx_{t+1}-\yy_t}^2\\
      &= f(\xx_{t+1}) + \beta_f(\xx_{t+1},\yy_t) -\frac{La_{t+1}^2}{A_{t+1}^2}\norm{\vv_{t+1}-\vv_t}^2\\
      &\geq f(\xx_{t+1}) + \beta_f(\xx_{t+1},\yy_t) -\frac{2La_{t+1}^2}{A_{t+1}^2}\norm{\wtilde\vv_{t+1}-\wtilde \vv_t}^2 - \frac{2La_{t+1}^2}{A_{t+1}^2}\norm{\ee_{t+1}-\ee_t}^2
  \end{align*}
  We plug it back in:
  \begin{align*}
      & \Eb{A_tf(\xx_t) + a_{t+1}f(\xx^\star) + \frac{1}{2}\norm{\wtilde\vv_t-\xx^\star}^2}\\
  = & \Eb{A_{t+1}f(\xx_{t+1}) + A_{t+1}\beta_f(\xx_{t+1},\yy_t) + A_t\beta_f(\yy_t,\xx_t) + a_{t+1}\beta_f(\yy_t,\xx^\star)  }\\
  & + \Eb{(\frac{1}{2}-\frac{2La_{t+1}^2}{A_{t+1}})\norm{\wtilde\vv_t-\wtilde\vv_{t+1}}^2+ \frac{1}{2}\norm{\wtilde\vv_{t+1}-\xx}^2 + a_{t+1}\inp{\gg_t -\nabla f(\yy_t)}{\wtilde\vv_{t+1}-\wtilde\vv_t}  } \\
  & - \left(\Eb{a_{t+1}\inp{\nabla f(\yy_t)}{\ee_{t+1}} + \frac{2La_{t+1}^2}{A_{t+1}}\norm{\ee_{t+1}-\ee_t}^2}\right)\\
  \overset{(i)}{\geq} & \Eb{A_{t+1}f(\xx_{t+1}) + A_{t+1}\beta_f(\xx_{t+1},\yy_t) + A_t\beta_f(\yy_t,\xx_t) + a_{t+1}\beta_f(\yy_t,\xx^\star) }\\
  & + \Eb{(\frac{1}{2}-\frac{2La_{t+1}^2}{A_{t+1}})\norm{\wtilde\vv_t-\wtilde\vv_{t+1}}^2 + \frac{1}{2}\norm{\wtilde\vv_{t+1}-\xx}^2 } \\
  & - \left(\Eb{\frac{3a_{t+1}^2\sigma_{\gg}^2}{2} + \frac{1}{6}\norm{\wtilde\vv_{t+1}-\wtilde\vv_t}^2 +  a_{t+1}\inp{\nabla f(\yy_t)}{\ee_{t+1}}+ \frac{2La_{t+1}^2}{A_{t+1}}\norm{\ee_{t+1}-\ee_t}^2}\right)
  \end{align*}
  where in (i) we applied Young's inequality and the bounded variance assumption. Now we pick $\xx\eqdef\xx^\star$, and subtract $A_{t+1}f(\xx^\star)$ from both sides, and rearrange the terms. We get:
  \begin{align*}
      A_{t+1}F_{t+1} &\leq A_tF_t + \frac{1}{2}(\wtilde R_t^2-\wtilde R_{t+1}^2) + \frac{3a_{t+1}^2\sigma_{\gg}^2}{2}- \left(\frac{1}{3}-\frac{2La_{t+1}^2}{A_{t+1}}\right)\EbNSq{\wtilde\vv_t-\wtilde\vv_{t+1}} - a_{t+1}\Eb{\beta_f(\yy_t,\xx^\star)}\\
      &\quad -\Eb{A_t\beta_f(\yy_t,\xx_t)} - \Eb{A_{t+1}\beta_f(\xx_{t+1},\yy_t)} +\Eb{a_{t+1}\inp{\nabla f(\yy_t)}{\ee_{t+1}}+ \frac{2La_{t+1}^2}{A_{t+1}}\norm{\ee_{t+1}-\ee_t}^2}
  \end{align*}

  Now we need to upper bound $a_{t+1}\EbInp{\nabla f(\yy_t)}{\ee_{t+1}}$. There are two options:
  \begin{align*}
      a_{t+1}\EbInp{\nabla f(\yy_t)}{\ee_{t+1}} &= a_{t+1}\EbInp{\gg_t }{\ee_{t+1}} + a_{t+1}\EbInp{\nabla f(\yy_t)-\gg_t }{\ee_{t+1}}\\
      &\leq \frac{1}{6}\EbNSq{\wtilde\vv_t-\wtilde\vv_{t+1}} + \frac{a_{t+1}^2\sigma_{\gg}^2}{2} + 2E_{t+1}
  \end{align*}
  The other option is:
  \begin{align*}
      a_{t+1}\EbInp{\nabla f(\yy_t)}{\ee_{t+1}} &\leq \frac{a_{t+1}}{4L}\EbNSq{\nabla f(\yy_t)} + a_{t+1}LE_{t+1}\\
      &\leq \frac{a_{t+1}}{2}\beta_f(\yy_t,\xx^\star) + a_{t+1}LE_{t+1}
  \end{align*}
  We combine these two upper bounds and get:
  \[
      a_{t+1}\EbInp{\nabla f(\yy_t)}{\ee_{t+1}} \leq \frac{1}{6}\EbNSq{\wtilde\vv_t-\wtilde\vv_{t+1}} + \frac{a_{t+1}^2\sigma_{\gg}^2}{2}+\frac{a_{t+1}}{2}\beta_f(\yy_t,\xx^\star) + \min\{2,a_{t+1}L\}E_{t+1}
  \]
  Plug it back in, we get:
  \begin{align*}
      A_{t+1}F_{t+1} &\leq A_tF_t + \frac{1}{2}(\wtilde R_t^2-\wtilde R_{t+1}^2) +2a_{t+1}^2\sigma_{\gg}^2- \left(\frac{1}{6}-\frac{2La_{t+1}^2}{A_{t+1}}\right)\EbNSq{\wtilde\vv_t-\wtilde\vv_{t+1}} - \frac{a_{t+1}}{2}\Eb{\beta_f(\yy_t,\xx^\star)}\\
      &\quad -\Eb{A_t\beta_f(\yy_t,\xx_t)} - \Eb{A_{t+1}\beta_f(\xx_{t+1},\yy_t)} +\left(\min\{2,a_{t+1}L\} + \frac{4La_{t+1}^2}{A_{t+1}}\right)E_{t+1} + \frac{4La_{t+1}^2}{A_{t+1}}E_t
  \end{align*}

  Now we sum over both sides from $t=0$ to $T-1$, and noticing that $E_t=0$ and $R_0^2=\wtilde R_0^2$, we get the desired result.
\end{proof}

\section{Missing Proofs for \Cref{sec:acc-ef}}
\label{sec:missing-proofs}
In this section we present the missing proofs for the analysis of \Cref{alg:dist-acc-ef}.

\eDescentSumInH*

\begin{proof}
    We first prove \Cref{eq:error-definition} by induction. For the base case $t=0$, we have $\avg{i}{n}\ee_0^i=\ee_0=\0$. Now assume that the lemma holds for some $t\geq 0$, we have:
    \begin{align*}
        \sum_{j=0}^{t}a_{j+1}(\hat \gg_j-\gg_j) &= \sum_{j=0}^{t-1}a_{j+1}(\hat\gg_j-\gg_j) + a_{t+1}(\hat \gg_t - \gg_t )\\
        &\overset{(i)}{=} \avg{i}{n}\ee_t^i+ a_{t+1}(\hat \gg_t - \gg_t )\\
        &= a_{t+1}\avg{i}{n}( \Delta_t^i - \gg_t^i +\frac{1}{a_{t+1}}\ee_t^i +\wtilde\gg_t^i)\\
        &=\avg{i}{n}\ee_{t+1}^i,
    \end{align*}
    where in $(i)$ we used the induction hypothesis. This completes the proof for \Cref{eq:error-definition}. This identity allows us to loosen the upper bound on $E_t$, and we consider the following quantity:
    \begin{equation}
        \label{eq:bar-e}
        \bar E_t\eqdef \avg{i}{n}\EbNSq{\ee_t^i}
    \end{equation}
    We note that $\bar E_t$ is an upper bound on $E_t$, we have:
    \[
        E_t = \EbNSq{\avg{i}{n}\ee_t^i}\leq \avg{i}{n}\EbNSq{\ee_t^i} = \bar E_t
    \]
    Therefore, instead of upper bounding $E_t$ directly, we will consider upper bounding $\bar E_t$. To utilize \Cref{lem:sum-lemma} to upper bound the weighted some of $\bar E_t$, we first give an individual descent on $\bar E_t$:

    \begin{align*}
        \bar E_{t+1} &= a_{t+1}^2\avg{i}{n}\EbNSq{\Delta_t^i-\delta_t^i}\\
            &\overset{(ii)}{\leq} (1-\delta)a_{t+1}^2\avg{i}{n}\EbNSq{ \gg_t ^i - \wtilde\gg_t ^i-\frac{1}{a_{t+1}}\ee_t^i}\\
            &\overset{(iii)}{\leq} (1-\frac{\delta}{2})\bar E_t + \frac{2a_{t+1}^2}{\delta}\avg{i}{n}\EbNSq{\gg_t ^i - \wtilde\gg_t ^i}
    \end{align*}
    where in $(ii)$ we used the definition of $\cC$. In $(iii)$ we used the Young's inequality. Apply the definition of $H_t$ and we get:
    \begin{equation}
        \label{eq:e-descent-one-step}
        \bar E_{t+1} \leq (1-\frac{\delta}{2})\bar E_t + \frac{2a_{t+1}^2}{\delta }H_t
    \end{equation}
    For $t\in[T-1]$, we have:
    \begin{align*}
        w_{t+1} \bar E_t &\leq (1-\frac{\delta}{2})w_{t+1}\bar E_t + \frac{2w_{t+1}'}{\delta}H_t\\
            &\overset{(iv)}{\leq }(1-\frac{\delta}{4})w_t\bar E_t + \frac{2w_{t+1}'}{\delta}H_t
    \end{align*}
    where in $(iv)$ we used the assumption that $w_{t+1}'(1-\frac{\delta}{2}) \leq w_t'(1-\frac{\delta}{4})$. Now we note that $\bar E_0=E_0=0$, write $\alpha_t = w_t\bar E_t,\lambda_t= w_{t+1}'H_t, \gamma = \frac{\delta}{4}$ and apply \Cref{lem:sum-lemma} and get the desired result.
\end{proof}

\hDescentSum*
\begin{proof}
    We first give a one step descent on $H_t$:
    \begin{align*}
        H_{t+1} &= \avg{i}{n}\EbNSq{\gg_{t+1}^i - \wtilde\gg_{t+1}^i}\\
            &= \avg{i}{n}\EbNSq{\gg_{t+1}^i - \wtilde \gg_t ^i - \wtilde\Delta_t^i}\\
            &\overset{(i)}{\leq } \frac{(1+\alpha)}{n}\sum_{i=1}^{n}\EbNSq{\wtilde \delta_t^i - \wtilde\Delta_t^i} + \frac{(1+\alpha^{-1})}{n}\sum_{i=1}^{n}\EbNSq{\gg_{t+1}^i - \gg_t ^i},\,
    \end{align*}
    where in $(i)$ we used Young's inequality. Pick $\alpha = \frac{\delta}{2(1-\delta)}$, we have:
    \begin{align*}
        H_{t+1} &\overset{(ii)}{\leq} (1-\frac{\delta}{2})H_t + \frac{2}{\delta n}\sum_{i=1}^{n}\EbNSq{\gg_{t+1}^i - \gg_t ^i}\\
        &\overset{(iii)}{\leq} (1-\frac{\delta}{2})H_t + \frac{4}{\delta n}\sum_{i=1}^{n}\EbNSq{\nabla f_i(\yy_{t+1}) - \nabla f_i(\yy_t)} + \frac{4\sigma^2}{\delta},\,
    \end{align*}
    where in $(ii)$ we used the definition of $\cC$. In $(iii)$ we used the Young's Inequality and \Cref{assumption:bounded_variance}. 

    For all $t\in[T-2]$ we have:
    \begin{align*}
        w_{t+2}'H_{t+1}&\leq w_{t+2}'(1-\frac{\delta}{2})H_t + \frac{4 w_{t+2}'}{\delta n}\sum_{i=1}^{n}\EbNSq{\nabla f_i(\yy_{t+1}) - \nabla f_i(\yy_t)} + \frac{4w_{t+2}'\sigma^2}{\delta}\\
        &\overset{(iv)}{\leq}w_{t+1}'(1-\frac{\delta}{4})H_t + \frac{4w_{t+2}'}{\delta n}\sum_{i=1}^{n}\EbNSq{\nabla f_i(\yy_{t+1}) - \nabla f_i(\yy_t)} + \frac{4w_{t+2}'\sigma^2}{\delta},\,
    \end{align*}
    where $(iv)$ uses the assumption. Note that by the algorithm we have $H_0=0$. Now we take $w_{t+1}'H_t$ to be $\alpha_t$, $\frac{\delta}{4}$ to be $\gamma$, and $\frac{4 w_{t+2}'}{\delta n}\sum_{i=1}^{n}\EbNSq{\nabla f_i(\yy_{t+1}) - \nabla f_i(\yy_t)} + \frac{4w_{t+2}'\sigma^2}{\delta}$ to be $\lambda_t$ in \Cref{lem:sum-lemma}, we can then apply \Cref{lem:sum-lemma} to get the desired result.
\end{proof}

\convergenceRate*
\begin{proof}

    We first further upper bound the weighted sum of $\bar E_t$ by noticing the following;
    \begin{align*}
        \avg{i}{n}\EbNSq{\nabla f_i(\yy_{t+1}) - \nabla f_i(\yy_t)} & \overset{(i)}{\leq} \frac{2}{n}\sum_{i=1}^{n}\left(\EbNSq{\nabla f_i(\yy_{t+1}) - \nabla f_i(\xx_t)} + \EbNSq{\nabla f_i(\xx_t)-\nabla f_i(\yy_t)}\right)\\
        &\overset{(ii)}{\leq} \frac{4\ell}{n}\left(\beta_f(\yy_{t+1},\xx_t) + \beta_f(\xx_t,\yy_t)\right),\,
    \end{align*}
    where in $(i)$ we used the Young's inequality and in $(ii)$ we used \Cref{assumption:f-i-smoothness}. We now apply this to \Cref{eq:e-descent-sum-in-h} from \Cref{lem:e-descent-sum-in-h} and get:
    \[
        \sum_{t=0}^{T-1}w_{t+1}\bar E_{t+1}
             \leq \frac{512\ell}{\delta^4 }\sum_{t=1}^{T-1}w_{t+1}'(\beta_f(\yy_{t},\xx_{t}) + \beta_f(\xx_{t},\yy_{t-1}))+ \frac{128\sigma^2}{\delta^4}\sum_{t=1}^{T-1}w_{t+1}'
    \]

    Now consider $a_t \eqdef \frac{t+\frac{32}{\delta}}{M}$ and $A_0 = \frac{32^2}{2\delta M}$. It's easy to verify that $w_{t+1}'(1-\frac{\delta}{2})\leq w_t'(1-\frac{\delta}{4})$ and $w_{t+1}(1-\frac{\delta}{2})\leq w_t(1-\frac{\delta}{4})$.

    We combine the above with \Cref{eq:main-descent} from \Cref{thm:main-descent} and get the following:
    \begin{align*}
        A_TF_T &\leq A_0F_0+\frac{R_0^2}{2} + \frac{2a_T^2T\sigma^2}{ n }+  \frac{512\ell}{\delta^4 }\sum_{t=1}^{T-1}w_{t+1}'(\beta_f(\yy_{t},\xx_{t}) + \beta_f(\xx_{t},\yy_{t-1}))+ \frac{128\sigma^2}{\delta^4}\sum_{t=1}^{T-1}w_{t+1}'\\
            &\quad  -\sum_{t=0}^{T-1}\left(\frac{1}{6}-\frac{2La_{t+1}^2}{A_{t+1}}\right)a_{t+1}^2\EbNSq{\gg_t }-\sum_{t=0}^{T-1}\left( \Eb{\frac{a_{t+1}}{2}\beta_f(\yy_t,\xx^\star) + A_t\beta_f(\yy_t,\xx_t) +A_{t+1}\beta_f(\xx_{t+1},\yy_t)}\right) \\
    \end{align*}
    We need:
    \[
        \frac{2La_{t+1}^2}{A_{t+1}}\leq \frac{1}{6},\quad\land\quad \frac{512\ell w_{t+1}'}{\delta^4}\leq A_t,\quad \land\quad \frac{512\ell w_{t+2}'}{\delta^4}\leq A_{t+1}
    \]

    When $M\geq 24L$, we have $\frac{2La_{t+1}^2}{A_{t+1}}\leq \frac{1}{6}$. For such $M$, we have $w_t\leq \frac{8}{3}$. It now suffices to have $M\geq \frac{2^{13}\ell}{\delta^4}$ for all the above inequalities to be satisfied. Therefore:
    \begin{align*}
        A_TF_T &\leq A_0F_0+\frac{R_0^2}{2} + \frac{2a_T^2T\sigma^2}{ n }+  \frac{128\sigma^2}{\delta^4}\sum_{t=1}^{T-1}w_{t+1}'
    \end{align*}
    Therefore:
    \[
        F_T\leq \frac{2^{10}F_0}{\delta^2\left(T+\frac{32}{\delta}\right)^2}+\frac{MR_0^2}{\left(T+\frac{32}{\delta}\right)^2} + \frac{4T\sigma^2}{Mn}+\frac{2^{9}LT\left(T+\frac{32}{\delta}\right)\sigma^2}{\delta^4M^2}
    \]
    Now we pick $M=\max\left\{24L,\frac{2^{13}\ell}{\delta^4}, \left(\frac{4T(T+\frac{32}{\delta})^2\sigma^2}{R_0^2n}\right)^{\frac{1}{2}}, 8\left(\frac{LT(T+\frac{32}{\delta})^3\sigma^2}{\delta^4R_0^2}\right)^{\frac{1}{3}}\right\}$, we get:
    \[
        F_T\leq \frac{24LR_0^2}{(T+\frac{32}{\delta})^2}+ \frac{2^{10}F_0}{\delta^2\left(T+\frac{32}{\delta}\right)^2}+\frac{2^{13}\ell R_0^2}{\delta^4\left(T+\frac{32}{\delta}\right)^2} + \frac{8L^{\frac{1}{3}}\sigma^{\frac{2}{3}}R_0^{\frac{4}{3}}}{\delta^{\frac{4}{3}}T^{\frac{2}{3}}}+ \frac{2R_0\sigma}{\sqrt{T}n} 
    \]
    It implies that, to reach $\varepsilon$-suboptimality, i.e. $F_T\leq \epsilon$, we need:
    \[
        T = \cO\left(\frac{R_0^2\sigma^2}{n\varepsilon^2}+ \frac{\sqrt{L}R_0^2\sigma}{\delta^2\varepsilon^{\nicefrac{3}{2}}}+\frac{\sqrt{\ell R_0^2}}{\delta^2\sqrt{\varepsilon}} + \frac{\sqrt{LR_0^2}}{\delta \sqrt{\varepsilon}}\right),\,
    \]
    Note that $L\leq \ell$ by \Cref{fact:ell-L} and we get the desired result.
\end{proof}

\section{Convergence of Acceleration with Vanilla EF}
\label{sec:acceleration-vanilla-ef}

\begin{algorithm}[tb]
    \caption{Accelerated Distributed Vanilla Error Feedback}
    \label{alg:dist-acc-ef-vanilla}
    \begin{algorithmic}[1]
        \State \textbf{Input:} $\xx_0, \yy_0, \vv_0, \ee_0^i=\0, A_0, (a_t)_{t=1}^\infty$, and $\cC$
        \For{$t = 0,1,\dots$}
        \State {\bf server side:}
        \State $A_{t+1} = A_t+a_{t+1}$
        \State $\yy_t = \frac{A_t}{A_{t+1}}\xx_t+ \frac{a_{t+1}}{A_{t+1}}\vv_t$
        \State server send $\xx_t,\yy_t$ to the clients
        \State {\bf each client $i$:}
        \State $\hat \gg_t ^i = \cC(\frac{1}{a_{t+1}} \ee_t^i+ \gg_t ^i)$
        \State  $\ee_{t+1}^i = \ee_t^i + a_{t+1} ( \gg_t ^i - \hat \gg_t ^i)$
        \State send to server $\hat \gg_t ^i$ 
        \State {\bf server side:}
        \State $\hat \gg_t ^i=\avg{i}{n}\hat\gg_t ^i$
        \State $\vv_{t+1}= \vv_t -a_{t+1} \gg_t ^i$
        \State $\xx_{t+1} = \frac{A_t}{A_{t+1}}\xx_t + \frac{a_{t+1}}{A_{t+1}}\vv_{t+1}$
        \hfill 
        \EndFor
    \end{algorithmic}	
\end{algorithm}

In this section, we consider an algorithm that combines acceleration with vanilla error feedback, summarized in \Cref{alg:dist-acc-ef-vanilla}. We apply \Cref{thm:main-descent} and give a tailored upper bound on the weighted sum of $E_t$ (in fact, as before, we upper bound the sum of $\bar E_t$). We need the following gradient similarity assumption:

\begin{assumption}
    \label{assumption:gradient-similarity}
    We assume that there exists $\zeta^2$ such that for all $\xx\in \R^d$, it holds:
    \begin{equation}
        \label{eq:gradient-similarity}
        \avg{i}{n}\norm{\nabla f_i(\xx)-\nabla f(\xx)}^2 \leq \zeta^2
    \end{equation}
\end{assumption}

Since $(\ee_t^i)_{t=0}^\infty$ also follows \Cref{lem:e-descent-sum-in-h}, it suffices that we give the following simple upper bound on $H_t$ where $\wtilde\gg_t^i=\0$:

\begin{lemma}
    \label{lem:h-bound-vanilla}
    Given \Cref{assumption:convexity,assumption:smoothness,assumption:bounded_variance,assumption:gradient-similarity}, for all $t\geq 0$, it holds that:
    \begin{equation}
        \label{eq:h-bound-vanilla}
        H_t \leq 2\zeta^2 + 4L\Eb{\beta_f(\yy_t,\xx^\star)} + \sigma^2
    \end{equation}
\end{lemma}
\begin{proof}
    \begin{align*}
        H_t &= \avg{i}{n}\EbNSq{\gg_t^i}\\
        &\overset{(i)}{\leq} \avg{i}{n}\EbNSq{\nabla f_i(\yy_t)} + \sigma^2\\
        &\overset{(ii)}{\leq}\frac{2}{n}\sum_{i=1}{n}\EbNSq{\nabla f_i(\yy_t)-\nabla f(\yy_t)} + 2\EbNSq{\nabla f(\yy_t)} + \sigma^2\\
        &\overset{(iii)}{\leq}2\zeta^2 + 4L\Eb{\beta_f(\yy_t,\xx^\star)} + \sigma^2,\,
    \end{align*}
    where in $(i)$ we applied \Cref{assumption:bounded_variance}, in $(ii)$ we used the Young's inequality, and in $(iii)$ we used \Cref{assumption:convexity,assumption:gradient-similarity,assumption:smoothness}.
\end{proof}

Now to give an upper bound on $\sum_{t=0}^{T-1}w_{t+1}E_{t+1}$:

\begin{lemma}
    \label{lem:e-descent-sum-vanilla}
    Given \Cref{assumption:convexity,assumption:smoothness,assumption:bounded_variance,assumption:gradient-similarity}, for all $T\geq 1$, if $w_{t+1}(1-\frac{\delta}{2})\leq w_t(1-\frac{\delta}{4})$, it holds that:
    \begin{equation}
        \label{eq:e-descent-sum-vanilla}
        \sum_{t=0}^{T-1} w_{t+1} E_{t+1} \leq  \frac{32L}{\delta^2 } \sum_{t=0}^{T-1}  w_{t+1}'\Eb{\beta_f(\yy_t,\xx^\star)} + \frac{4(4\zeta^2+\delta\sigma^2)}{\delta^2} \sum_{t=0}^{T-1}  w_{t+1}'
    \end{equation}
\end{lemma}

\begin{proof}
    We simply combine \Cref{eq:e-descent-sum-in-h} from \Cref{lem:e-descent-sum-in-h} and \Cref{eq:h-bound-vanilla} from \Cref{lem:h-bound-vanilla} and get the desired result.

\end{proof}

Now we can combine \Cref{eq:e-descent-sum-vanilla} from \Cref{lem:e-descent-sum-vanilla} with \Cref{eq:main-descent} from \Cref{thm:main-descent} to get the following convergence rate:

\begin{theorem}
    \label{thm:convergence-rate-vanilla}
    Given \Cref{assumption:convexity,assumption:smoothness,assumption:bounded_variance,assumption:gradient-similarity}, and let $a_t\eqdef \frac{t+\frac{4}{t}}{M}$ and $A_0\eqdef \frac{8}{\delta^2 M}$, it suffices to have:
    \begin{equation}
        \label{eq:convergence-rate-vanilla}
        T = \cO\left(\frac{\sigma^2 R_0^2}{\varepsilon^2} + \frac{\sqrt{L}R_0^2(\nicefrac{\zeta}{\sqrt{\delta}}+\sigma)}{\sqrt{\delta}\varepsilon^{\nicefrac{3}{2}}}+ \frac{LR_0^2}{\delta\varepsilon }\right)
    \end{equation}
    number of iterations of \Cref{alg:dist-acc-ef-vanilla} to get $F_T\leq \varepsilon$, where we can set 
    \[
        M=\max\left\{\frac{40L(T+\frac{4}{\delta})}{\delta},\sqrt{\frac{4T(T+\frac{4}{\delta})^2\sigma^2}{R_0^2n}}, \left(\frac{544L(4\zeta^2+\delta\sigma^2)T(T+\frac{4}{\delta})^3}{\delta^2R_0^2}\right)^{\frac{1}{3}}\right\}.
    \]
\end{theorem}

\begin{proof}
    It's easy to check that for this choice of $a_t$ and $A_0$, we have that $w_{t+1}(1-\frac{\delta}{2})\leq w_t(1-\frac{\delta}{4})$.

    In addition, we have that $\frac{La_{t+1}^2}{A_{t+1}}\leq \frac{2L}{M}$.

    We combine \Cref{eq:main-descent} from \Cref{thm:main-descent} and \Cref{eq:e-descent-sum-vanilla} from \Cref{lem:e-descent-sum-vanilla}, pick $M\geq 24L$ and get the following:
    \begin{align*}
        A_TF_T 
            &\leq A_0F_0+\frac{R_0^2}{2} + \frac{2a_T^2T\sigma^2}{ n }\\
            &\quad + \frac{32L}{\delta^2 } \sum_{t=0}^{T-1}  w_{t+1}'\Eb{\beta_f(\yy_t,\xx^\star)} + \frac{8(2\zeta^2+\sigma^2)}{\delta^2 } \sum_{t=0}^{T-1}  w_{t+1}'\\
            &\quad -\sum_{t=0}^{T-1}\left( \Eb{\frac{a_{t+1}}{2}\beta_f(\yy_t,\xx^\star) + A_t\beta_f(\yy_t,\xx_t) + A_{t+1}\beta_f(\xx_{t+1},\yy_t)}\right)  \\
            &\leq A_0F_0+\frac{R_0^2}{2} + \frac{2a_T^2T\sigma^2}{ n }\\
            &\quad + \frac{32L^2}{\delta^2 M^3} \sum_{t=0}^{T-1}(t+17+\frac{4}{\delta})(t+1+\frac{4}{\delta})^2\Eb{\beta_f(\yy_t,\xx^\star)}\\
            & \quad + \frac{4L(4\zeta^2+\delta\sigma^2)}{\delta^2 M^3} \sum_{t=0}^{T-1}(t+17+\frac{4}{\delta})(t+1+\frac{4}{\delta})^2\\
            &\quad -\sum_{t=0}^{T-1}\left( \Eb{\frac{a_{t+1}}{2}\beta_f(\yy_t,\xx^\star) + A_t\beta_f(\yy_t,\xx_t) + A_{t+1}\beta_f(\xx_{t+1},\yy_t)}\right) 
    \end{align*}
    We need $M$ such that:
    \[  
        \frac{32L^2(t+17+\frac{4}{\delta})(t+1+\frac{4}{\delta})^2}{\delta^2 M^2}\leq \frac{t+1+\frac{4}{\delta}}{2}, \forall 0\leq t\leq T-1
    \]  
    It suffices that $M\geq \frac{40L(T+\frac{4}{\delta})}{\delta}$, and we have the following:
    \[
        F_T \leq \frac{16F_0}{\delta^2 (T+\frac{4}{\delta})^2}+\frac{MR_0^2}{(T+\frac{4}{\delta})^2} + \frac{4T\sigma^2}{Mn} + \frac{544L(4\zeta^2+\delta\sigma^2)T(T+\frac{4}{\delta})}{\delta^2M^2}
    \]
    Now pick $M=\max\{\frac{40L(T+\frac{4}{\delta})}{\delta},\sqrt{\frac{4T(T+\frac{4}{\delta})^2\sigma^2}{R_0^2n}}, \left(\frac{544L(4\zeta^2+\delta\sigma^2)T(T+\frac{4}{\delta})^3}{\delta^2R_0^2}\right)^{\frac{1}{3}}\}$, we would have that after 
    \[
        T = \cO\left(\frac{\sigma^2 R_0^2}{\varepsilon^2} + \frac{\sqrt{L}R_0^2(\nicefrac{\zeta}{\sqrt{\delta}}+\sigma)}{\sqrt{\delta}\varepsilon^{\nicefrac{3}{2}}}+ \frac{LR_0^2}{\delta\varepsilon }\right)
    \]
    iterations of \Cref{alg:dist-acc-ef-vanilla}, we have $F_T\leq \varepsilon$.
\end{proof}

\begin{remark}
    \label{remark:convergence-rate-vanilla}
    We notice that \Cref{alg:dist-acc-ef-vanilla} fail to achieve accelerated rate in the determinstic term and obtain the same rate as the unaccelerated error feedback. This is due to the lower bound where we need $M\geq \frac{40L(T+\frac{4}{\delta})}{\delta}$, which canceled out the acceleration effect. Such a requirement on the minimal $M$ comes the fact that in \Cref{eq:main-descent}, we only have $\frac{a_{t+1}}{2}\Eb{\beta_f(\yy_t,\xx^\star)}$ to cancel out the error terms involving $\beta_f(\yy_t,\xx^\star)$. In contrast, in \Cref{alg:dist-acc-ef}, we applied the gradient difference compression, which enables us to the cancel out the errors with $A_t\Eb{\beta_f(\yy_t,\xx_t)}$ and $A_{t+1}\Eb{\beta_f(\xx_{t+1},\yy_t)}$ and therefore achieve the accelerated $\cO(\frac{1}{\sqrt{\varepsilon}})$ rate.
\end{remark}

\section{Convergence of Accelerated Distributed Optimization with Absolute Compression}
\label{sec:absolute-compression}

\begin{algorithm}[tb]
    \caption{Accelerated Distributed Optimization with Absolute Compression}
    \label{alg:dist-acc-absolute}
    \begin{algorithmic}[1]
        \State \textbf{Input:} $\xx_0, \yy_0, \vv_0, \ee_0^i=\0, A_0, (a_t)_{t=1}^\infty$, and $\cC_{\Delta}$
        \For{$t = 0,1,\dots$}
        \State {\bf server side:}
        \State $A_{t+1} = A_t+a_{t+1}$
        \State $\yy_t = \frac{A_t}{A_{t+1}}\xx_t+ \frac{a_{t+1}}{A_{t+1}}\vv_t$
        \State server send $\xx_t,\yy_t$ to the clients
        \State {\bf each client $i$:}
        \State $\hat \gg_t ^i = \cC_{\Delta}(\gg_t ^i)$
        \State send to server $\hat \gg_t ^i$ 
        \State {\bf server side:}
        \State $\hat \gg_t ^i=\avg{i}{n}\hat\gg_t ^i$
        \State $\vv_{t+1}= \vv_t -a_{t+1} \gg_t ^i$
        \State $\xx_{t+1} = \frac{A_t}{A_{t+1}}\xx_t + \frac{a_{t+1}}{A_{t+1}}\vv_{t+1}$
        \hfill 
        \EndFor
    \end{algorithmic}	
\end{algorithm}

In this section we consider a different kind of compressor, absolute compression:

\begin{definition}
    \label{def:absolute-compression}
    We say that a (possibly randomized) mapping  $\cC_{\Delta} \colon \R^d \to \R^d$ is an $\Delta$-absolute compression operator if for some constant $\Delta\geq 0$ it holds 
    \begin{equation}
        \label{eq:absolute-compression}
        \Eb{\norm{\cC(\xx) -\xx}^2} \leq \Delta^2\, \quad \forall \xx \in \R^d. 
    \end{equation}
\end{definition}
This class of compressor also include some popular examples of compressions, including hard-thresholding~\citep{sahu2021rethinking}, rounding with bounded errors~\citep{gupta2015deep}, and scaled integer rounding~\citep{sapio2021scaling}.

We summarize the accelerated distributed optimization with absolute compression method in \Cref{alg:dist-acc-absolute}. Here the algorithm simply aggregates the compressed gradients and perform the update. We demonstrate our framework with \Cref{thm:main-descent} is very flexible and can also be applied to analyze \Cref{alg:dist-acc-absolute}.

We pick the standard choice of $a_t$ and $A_0$: $a_t\eqdef \frac{t}{M}$ and $A_0=0$. We can apply \Cref{thm:main-descent} to get the following convergence rate:
\begin{theorem}
    \label{thm:absolute-compression-rate}
    Given \Cref{assumption:convexity,assumption:smoothness,assumption:bounded_variance}, and let $a_t\eqdef \frac{t}{M}$ and $A_0=0$, then we have
    \begin{equation}
        \label{eq:absolute-compression-rate}
        F_T\leq \frac{24LR_0^2}{T^2} + \sqrt{\frac{4R_0^2\sigma^2}{nT} }+ \left(32LR_0^4\Delta^2\right)^{\frac{1}{3}},\,
    \end{equation}
    where we set $M=\max\{24L, \sqrt{\frac{4T^3\sigma^2}{nR_0^2}}, \left(\frac{2L\Delta^2(T+16)T^5}{R_0^2}\right)^{\frac{1}{3}}\}$.
\end{theorem}
\begin{proof}
    With the choice of $a_t$, we can give an upper bound on $E_t$:
    \begin{align*}
        E_{t} & = \EbNSq{\sum_{j=0}^{t-1}a_{j+1}(\hat \gg_j -\gg_j )}\\
            & = \EbNSq{\sum_{j=0}^{t-1}a_{j+1}(\cC_{\Delta}(\gg_j)-\gg_j)}\\
            &\leq a_{t+1}^2\EbNSq{\sum_{j=0}^{t-1}(\cC_{\Delta}(\gg_j)-\gg_j)}\\
            &\leq a_{t+1}^2t\sum_{j=0}^{t-1}\EbNSq{\cC_{\Delta}(\gg_j)-\gg_j}\\
            &\leq a_{t+1}^2t^2\Delta^2 
    \end{align*}
    Therefore:
    \begin{align*}
        \sum_{i=0}^{T-1}w_{t+1}E_{t+1} &\leq \frac{L\Delta^2}{M^3} \sum_{i=0}^{T-1} (t+17)(t+1)^2t^2\\
        &\leq \frac{L\Delta^2(T+16)T^5}{M^3}
    \end{align*}
    Now we plug it into \Cref{eq:main-descent} from \Cref{thm:main-descent}, assume that $M\geq 24L$, we get the following:
    \begin{align*}
        A_TF_T &\leq A_0F_0+\frac{R_0^2}{2} + \frac{2a_T^2T\sigma^2}{ n } + \frac{L\Delta^2(T+16)T^5}{M^3}\\
        &= \frac{R_0^2}{2} + \frac{2T^3\sigma^2}{M^2n} + \frac{L\Delta^2(T+16)T^5}{M^3}
    \end{align*}
    Dividing both sides by $A_T$, we get:
    \[
        F_T\leq \frac{MR_0^2}{T^2} + \frac{4T\sigma^2}{Mn} + \frac{2L\Delta^2(T+16)T^3}{M^2}
    \]
    Now pick $M=\max\{24L, \sqrt{\frac{4T^3\sigma^2}{nR_0^2}}, \left(\frac{2L\Delta^2(T+16)T^5}{R_0^2}\right)^{\frac{1}{3}}\}$, we would have that:
    \[
        F_T\leq \frac{24LR_0^2}{T^2} + \sqrt{\frac{4R_0^2\sigma^2}{nT} }+ \left(32LR_0^4\Delta^2\right)^{\frac{1}{3}}.
    \]
\end{proof}
\begin{remark}
    We see that \Cref{alg:dist-acc-absolute} converges up to some $\cO(LR_0^4\Delta^2)$ neighborhood of the optimal solution. There exists some works that scales the compressed message by the inverse of the stepsize, see, e.g. \citet{fatkhullin2023momentum}. We do not pursue this direction as here we are interested in the applicability of our framework \Cref{thm:main-descent}, instead of the best possible rate.
\end{remark}

\section{Method of Repeated Communication}
\label{sec:repeated-communication}

\begin{figure}[tb]
    \begin{minipage}{0.49\textwidth}
        \begin{figure}[H]
            \begin{algorithm}[H]
                \caption{Repeated Compressor $\cC_R$}
                \label{alg:CR}
                \begin{algorithmic}[1]
                    \State \textbf{Input:} $\xx, \cC$, and $R$
                    \State \textbf{Client side:}
                    \State $\cc_0=\0$
                    \For{$q=1,2,\cdots,R$}
                        \State $\Delta_q = \cC(\xx-\cc_{q-1})$
                        \State $\cc_{q}=\cc_{q-1}+\Delta_q$
                    \hfill 
                    \EndFor
                    \State Client sends $\{\Delta_q\}_{q\in [R]}$ to server
                    \State \textbf{Server side:}
                    \State $\Delta = \sum_{q=1}^{R}\Delta_q$
                \end{algorithmic}	
            \end{algorithm}
        \end{figure}
    \end{minipage}
    \hfill
    \begin{minipage}{0.49\textwidth}
        \begin{figure}[H]
            \begin{algorithm}[H]
                \caption{\algname{NEOLITHIC}}
                \label{alg:neolithic}
                \begin{algorithmic}[1]
                    \State \textbf{Input:} $\xx_0,\yy_0,\vv_0, A_0, (a_t)_{t=1}^\infty, \cC_R$ and $M$
                    \For{$t = 0,1,\dots$}
                    \State {\bf server side:}
                    \State $A_{t+1} = A_t+a_{t+1}$
                    \State $\yy_t = \frac{A_t}{A_{t+1}}\xx_t+ \frac{a_{t+1}}{A_{t+1}}\vv_t$
                    \State {\bf each client $i$:}
                    \State $\gg_t ^i = \gg_i(\yy_t,\xi_{\yy_t}^i)$
                    \State $\hat\gg_t ^i = \cC_R(\gg_t ^i)$ 
                    \State send to server $\Delta_t^i$ 
                    \State {\bf server side:}
                    \State $\hat \gg_t =\avg{i}{n}\gg_t ^i$
                    \State $\vv_{t+1}= \vv_t -a_{t+1}\hat \gg_t $
                    \State $\xx_{t+1} = \frac{A_t}{A_{t+1}}\xx_t + \frac{a_{t+1}}{A_{t+1}}\vv_{t+1}$
                    \hfill 
                    \EndFor
                \end{algorithmic}	
            \end{algorithm}
        \end{figure}
    \end{minipage}
\end{figure}

In this section we further discuss \algname{NEOLITHIC} and the method of repeated communication, proposed by~\citet{huang2022lowerbounds} initially for unaccelerated methods, and later extended to accelerated methods by~\citet{he2023lower}. The method at its core is the same as compressed \algname{SGD} and compressed accelerated \algname{SGD}, where the server naively aggregates the compressed gradients and make an update. The only difference is that the compressor used here is a repeated compressor $\cC_R$ that compresses the gradient by applying the basic contractive compressor $\cC$ for $R$ rounds (and therefore communicating with the server for $R$ rounds). For completeness, we provide the precise definition of $\cC_R$ in \Cref{alg:CR}, which is taken (with some simple reformulation for clarity) from~\citet{he2023lower}.

The key point of the repeated compressor $\cC_R$ is that it reduces the error of the compression, we state the following lemma from~\citet{he2023lower}:
\begin{lemma}
    \label{lem:repeated-compressor}
    Given a $\delta$-contractive compressor $\cC$ and $R$, the repeated compressor $\cC_R$ satisfies the following:
    \begin{equation}
        \label{eq:repeated-compressor}
        \EbNSq{\cC_R(\xx)-\xx} \leq (1-\delta)^R\norm{\xx}^2,\forall \xx\in \R^d
    \end{equation}
\end{lemma}
However, we point out that when $R$ is some constant value independent of $\delta$ and $\delta$ is small enough (so that $R\leq\frac{1}{2\delta}$)
\[
    1-R\delta \leq (1-\delta)^R\leq (1-R\delta + R^2\delta^2)\leq (1-\frac{R\delta}{2})\,,
\]
implying that in such a case $\cC_R$ is essentially an $R\delta$-contractive compressor. Therefore, in such a case when $R$ is some fixed constant, the \algname{NEOLITHIC} algorithm simply reduces to the compressed \algname{SGD} with $R\delta$-contractive compressor. It is well known that compressed \algname{SGD} diverges~\citep{beznosikov2023biased} and therefore so will \algname{NEOLITHIC}. This also explains why in the practical situation where $R$ is some fixed constant, \algname{NEOLITHIC}'s performance is very poor, as observed in~\citep{fatkhullin2023momentum}. See also \Cref{sec:mnist-classification}

\citet{he2023lower} suggested choosing $R=\max\{\frac{4}{\delta}\ln(\frac{4}{\delta}), \frac{1}{\delta}\ln(24\kappa + \frac{25n^2\kappa^3\zeta^4}{\sigma^4}+5n\kappa^{\nicefrac{3}{2}})\}$ (where $\kappa$ denotes the total number of updates at the server), i.e. in between each updates, the algorithm sends $\max\{\frac{4}{\delta}\ln(\frac{4}{\delta}), \frac{1}{\delta}\ln(24\kappa + \frac{25n^2\kappa^3\zeta^4}{\sigma^4}+5n\kappa^{\nicefrac{3}{2}})\}$ messages compressed by $\cC$ to approximate the original gradient. Consider the classical example of Top-$K$ compressor, which is $\frac{K}{d}$-contractive compression, which sends $K$ entries out of $d$ entries of the uncompressed vector. In this case, $\max\{\frac{4}{\delta}\ln(\frac{4}{\delta}), \frac{1}{\delta}\ln(24\kappa + \frac{25n^2\kappa^3\zeta^4}{\sigma^4}+5n\kappa^{\nicefrac{3}{2}})\}$ Top-$K$ compressed messages will contain $\max\{4d\ln(\frac{4d}{K}),d\ln(\kappa)\}$ entries in total, much more than the original uncompressed message. Note also that running \Cref{alg:CR} with Top-$K$ for exactly $\frac{d}{K}$ rounds would send the uncompressed vector \emph{exactly}. More generally, if a $\delta$-contractive compression $\cC$ sends a $\tau$ fraction of bits of the uncompresed vector, then as shown in~\citet{zhu2023optimal}, in the worst case, 
\[
    \tau =\Omega(\delta)
\]
In other words, in the worst case, for any $\delta$-contractive compression, \algname{NEOLITHIC} sends a $\Omega(\max\{\ln(\frac{1}{\delta}),\ln(\kappa)\})$ factor of the size of the uncompressed local gradient. For reasonable $\delta$ and $\kappa$, this factor is larger than one, implying, again, there is no compression at all.

Finally, with such high accuracy compression (with the relative error $(1-\delta)^R\leq \frac{1}{\kappa}$), \citet{he2023lower} showed that the compressed accelerated (and the unaccelerated) \algname{SGD}'s convergence is essentialy unaffected by the compression error. Therefore, \algname{NEOLITHIC} performs, up to some constant, the same number of server side updates as the uncompressed method, that is, $\kappa = \cO(\frac{R_0^2\sigma^2}{n\varepsilon^2} + \frac{\sqrt{LR_0^2}}{\sqrt{\varepsilon}})$. Therefore, the total number of communication rounds becomes
\[
    R\cO\left(\frac{R_0^2\sigma^2}{n\varepsilon^2} + \frac{\sqrt{LR_0^2}}{\sqrt{\varepsilon}}\right) = \wtilde\cO\left(\frac{R_0^2\sigma^2}{\delta n\varepsilon^2} + \frac{\sqrt{LR_0^2}}{\delta\sqrt{\varepsilon}}\right)\,,
\]
where, importantly, the leading term has a $\wtilde\cO(\frac{1}{\delta})$ factor. Thererfore, \algname{NEOLITHIC} has no benefit over the uncompressed method (which is unsurprising as \algname{NEOLITHIC} itself is an uncompressed method in disguise). This is in stark contrast to most other works with Error Feedback, including ours, where the leading term is $\cO(\frac{R_0^2\sigma^2}{n\varepsilon^2})$. We point out that in the write up by~\citet{he2023lower}, they claimed that \algname{NEOLITHIC} achieves a $\wtilde\cO(\frac{R_0^2\sigma^2}{n\varepsilon^2})$ leading term. This is misleading, because, to achieve this rate, they have to take a $B=R=\max\{\frac{4}{\delta}\ln(\frac{4}{\delta}), \frac{1}{\delta}\ln(24\kappa + \frac{25n^2\kappa^3\zeta^4}{\sigma^4}+5n\kappa^{\nicefrac{3}{2}})\}$ batch for each update. Recall that taking a batch of size $B$ effectively divides the variance by $B$, and therefore, for a fair comparison, we point out that, our method (along with most other error feedback variants), would achieve a $o\left(\frac{\delta R_0^2\sigma^2}{n\varepsilon^2}\right)$ leading term when using the same batch size, still better by an $\Omega(\frac{1}{\delta})$ factor.

\subsection{Convergence Analysis of \algname{NEOLITHIC}}
In this section we provide an alternative analysis of \algname{NEOLITHIC} based on \Cref{thm:main-descent}, further demonstrating the flexibility of our proposed framework. \algname{NEOLITHIC} is summarized in \Cref{alg:neolithic}. We point out that the presentation of the algorithm here is slightly different from the original presentation in~\citet{he2023lower}, where we made some modifications to their underlying acceleration method to aligh with ours. Note that we also give slightly different recomendation on $R$, the number of rounds of communications between each server update.

Now we upper bound the weighted sum of $E_t$:
\begin{lemma}
    \label{lem:e-descent-sum-neolithic}
    Given \Cref{assumption:convexity,assumption:smoothness,assumption:bounded_variance,assumption:gradient-similarity}, we have:
    \begin{equation}
        \label{eq:e-descent-sum-neolithic}
        \begin{aligned}
            \sum_{t=0}^{T-1}w_{t+1}E_{t+1} &\leq (1-\delta)^R\zeta^2\sum_{t=0}^{T-1}w_{t+1}t\sum_{j=0}^{t-1}a_{j+1}^2 + 2L(1-\delta)^R\sum_{t=0}^{T-1}w_{t+1}t\sum_{j=0}^{t-1}a_{j+1}^2\beta_f(\yy_j,\xx^\star)\\
            &\quad + (1-\delta)^R\sigma^2\sum_{t=0}^{T-1}w_{t+1}t\sum_{j=0}^{t-1}a_{j+1}^2
        \end{aligned}
    \end{equation}
\end{lemma}
\begin{proof}
    \begin{align*}
        E_{t+1} &= \EbNSq{\sum_{j=0}^{t-1}a_{j+1} \gg_j-a_{j+1} \avg{i}{n}\hat\gg_j^i}\\
        &\leq t\sum_{j=0}^{t-1}a_{j+1}^2\EbNSq{\gg_j-\avg{i}{n}\hat\gg_j^i}\\
        &\leq \frac{t}{n}\sum_{j=0}^{t-1}a_{j+1}^2\sum_{i=1}^{n}\EbNSq{\gg_j-\hat\gg_t ^i}\\
        &\leq \frac{(1-\delta)^Rt}{n }\sum_{j=0}^{t-1}a_{j+1}^2\sum_{i=1}^{n}\EbNSq{\gg_j}\\
        &\leq \frac{(1-\delta)^Rt}{n }\sum_{j=0}^{t-1}a_{j+1}^2\sum_{i=1}^{n}\EbNSq{\nabla f_i(\yy_j)} + (1-\delta)^Rt\sigma^2\sum_{j=0}^{t-1}a_{j+1}^2\\
    \end{align*}
    We note the following:
    \begin{align*}
        \avg{i}{n}\EbNSq{\nabla f_i(\yy_j)} &= \avg{i}{n}\EbNSq{\nabla f_i(\yy_j)-\nabla f(\yy_j)} + \EbNSq{\nabla f(\yy_j)}\\
        &\leq \zeta^2 + 2L\beta_f(\yy_j,\xx^\star)
    \end{align*}
    Plugging it back in, we get:
    \[
        E_{t+1} \leq (1-\delta)^R\zeta^2t\sum_{j=0}^{t-1}a_{j+1}^2+2L(1-\delta)^Rt\sum_{j=0}^{t-1}a_{j+1}^2\beta_f(\yy_j,\xx^\star) + (1-\delta)^Rt\sigma^2\sum_{j=0}^{t-1}a_{j+1}^2
    \]
    Summing this up, we get the desired result.
\end{proof}

Now we can combine the \Cref{thm:main-descent} and \Cref{lem:e-descent-sum-neolithic} to get the convergence rate of \algname{NEOLITHIC}:

\begin{theorem}
    \label{thm:convergence-rate-neolithic}
    Given \Cref{assumption:convexity,assumption:smoothness,assumption:bounded_variance,assumption:gradient-similarity}, and let $a_t\eqdef \frac{t}{M},A_0=0$, and $R\eqdef \max\{\frac{4}{\delta}\ln(T), \frac{1}{\delta}\ln(\frac{4nT^2}{3}), \frac{1}{\delta}\ln(\frac{4n\zeta^2T^2}{3\sigma^2})\}$, it suffices to take 
    \begin{equation}
        \label{eq:convergence-rate-neolithic}
        RT = \cO\left(\left(\frac{R_0^2\sigma^2}{\delta n\varepsilon^2}+\frac{1}{\delta}\sqrt{\frac{LR_0^2}{ \varepsilon}}\right)\left(\ln(nT) + \ln(\frac{n\zeta^2 T^2}{\sigma^2})\right)\right)
    \end{equation}
    total number of communication rounds to get $F_T\leq \varepsilon$, where we can set $M=\max\{24L, \sqrt{\frac{12T^3\sigma^2}{nR_0^2}}\}$.
\end{theorem}

\begin{proof}
    As before, we assume that $M\geq 24L$. Now we have $w_t \leq  \frac{8}{3}$.
    \begin{align*}
        A_TF_T &\leq A_0F_0+\frac{R_0^2}{2} + \frac{2a_T^2T\sigma^2}{n }\\
            &\quad + (1-\delta)^R\zeta^2\sum_{t=0}^{T-1}w_{t+1}t\sum_{j=0}^{t-1}a_{j+1}^2 + 2L(1-\delta)^R\sum_{t=0}^{T-1}w_{t+1}t\sum_{j=0}^{t-1}a_{j+1}^2\beta_f(\yy_j,\xx^\star)\\
            &\quad + (1-\delta)^R\sigma^2\sum_{t=0}^{T-1}w_{t+1}t\sum_{j=0}^{t-1}a_{j+1}^2\\
            &\quad -\sum_{t=0}^{T-1}\left( \Eb{\frac{a_{t+1}}{2}\beta_f(\yy_t,\xx^\star) + A_t\beta_f(\yy_t,\xx_t) +A_{t+1}\beta_f(\xx_{t+1},\yy_t)}\right) \\
            &\leq A_0F_0+\frac{R_0^2}{2} + \frac{2a_T^2T\sigma^2}{ n }\\
            &\quad +  \frac{8(1-\delta)^R\zeta^2}{3M^2}\sum_{t=0}^{T-1}t^4 + \frac{16L(1-\delta)^RT^4}{3M^2}\sum_{t=0}^{T-1}\beta_f(\yy_t,\xx^\star)+ \frac{8(1-\delta)^R\sigma^2}{3M^2}\sum_{t=0}^{T-1}t^4\\
            &\quad -\sum_{t=0}^{T-1}\left( \Eb{\frac{a_{t+1}}{2}\beta_f(\yy_t,\xx^\star) + A_t\beta_f(\yy_t,\xx_t) +A_{t+1}\beta_f(\xx_{t+1},\yy_t)}\right) \\
            &\leq A_0F_0+\frac{R_0^2}{2} + \frac{2a_T^2T\sigma^2}{n }\\
            &\quad +  \frac{8(1-\delta)^R\zeta^2T^5}{3M^2}+ \frac{16L(1-\delta)^RT^4}{3M^2}\sum_{t=0}^{T-1}\beta_f(\yy_t,\xx^\star)+ \frac{8(1-\delta)^R\sigma^2T^5}{3M^2}\\
            &\quad -\sum_{t=0}^{T-1}\left( \Eb{\frac{a_{t+1}}{2}\beta_f(\yy_t,\xx^\star) + A_t\beta_f(\yy_t,\xx_t) +A_{t+1}\beta_f(\xx_{t+1},\yy_t)}\right) \\
    \end{align*}
    Now we need the following:
    \[
        \frac{16L(1-\delta)^RT^4}{3M} \leq \frac{t+1}{2}\quad \land\quad \frac{8(1-\delta)^RT^5}{3}\leq \frac{2T^3}{n}\quad \land \quad \frac{8(1-\delta)^R\zeta^2T^5}{3} \leq \frac{2T^3\sigma^2}{n}
    \]
    it suffices if $R\eqdef \max\{\frac{4}{\delta}\ln(T), \frac{1}{\delta}\ln(\frac{4nT^2}{3}), \frac{1}{\delta}\ln(\frac{4n\zeta^2T^2}{3\sigma^2})\}$. With this $R$, we have:
    \[
        A_TF_T\leq \frac{R_0^2}{2} + \frac{6T^3\sigma^2}{M^2n}
    \]
    Therefore, with $M=\max\{24L, \sqrt{\frac{12T^3\sigma^2}{nR_0^2}}\}$, we have that it suffices to take
    \[
        T = \cO(\frac{R_0^2\sigma^2}{n\varepsilon^2}+\sqrt{\frac{LR_0^2}{\varepsilon}})
    \]  
    number of server updates to achieve $\varepsilon$-accuracy, as desired.
\end{proof}

\end{document}